\pdfoutput=1

\documentclass[12pt]{amsart}
\usepackage[centertags]{amsmath}
\usepackage{amsfonts}
\usepackage{amssymb}
\usepackage{amsthm}
\usepackage{graphicx}
\usepackage{fullpage }

\usepackage[ small,nohug,heads=littlevee]{diagrams}

\diagramstyle[labelstyle=\scriptstyle]

\newtheorem{theorem}{Theorem}[section]
\newtheorem{proposition}[theorem]{Proposition}
\newtheorem{lemma}[theorem]{Lemma}
\theoremstyle{definition}

\theoremstyle{remark}
\newtheorem{remark}[theorem]{Remark}

\newcommand{\Z}{{\mathbb{Z}}}
\newcommand{\C}{{\mathcal{C}}}

\newcommand{\s}{\mathcal S}
\newcommand{\tC}{{\widetilde{\mathcal{C}}}}
\newcommand{\D}{\mathcal D}
\newcommand{\tD}{{\widetilde{\mathcal{D}}}}
\newcommand{\td}{{\widetilde{d}}}
\newcommand{\hD}{{\widehat{\mathcal{D}}}}
\newcommand{\hd}{{\widehat{d}}}
\newcommand{\hC}{{\widehat{\mathcal{C}}}}

\newcommand{\tH}{{\widetilde{H}}}
\newcommand{\hH}{\widehat{H}}

\newcommand{\hHB}{\widehat{HB}}
\newcommand{\tHK}{{\widetilde{HK}}}

\newcommand{\qdim}{\operatorname{qdim}}
\newcommand{\rk}{\operatorname{rk}}
\newcommand{\im}{\operatorname{im}}

\begin{document}
\title[The chromatic polynomial of fatgraphs and its categorification]
{The chromatic polynomial of fatgraphs and its categorification}

\author{Martin Loebl}
\address{Dept.~of Applied Mathematics and
Institute of Theoretical Computer Science \\
Charles University \\
Malostransk\' e n. 25 \\
118 00 Praha 1 \\ Czech Republic.
\newline  {\em and}  Depto. Ing. Matem\,atica \\ University of Chile \\ Chile.}

\email{loebl@kam.mff.cuni.cz}

\author{ Iain Moffatt}
\address{Dept.~of Applied Mathematics and
Institute of Theoretical Computer Science \\
Charles University \\
Malostransk\' e n. 25 \\
118 00 Praha 1 \\
Czech Republic.}
\curraddr{Department of Combinatorics and Optimization \\ University of Waterloo \\ Waterloo \\ Canada.}
\email{imoffatt@math.uwaterloo.ca}

\thanks{
\newline
2000 {\em Mathematics Classification.} Primary 05C10. Secondary 57M15.
}

\date{
This edition: November 26, 2007 \hspace{0.5cm} First edition: November 22, 2005.}

\begin{abstract}
Motivated by Khovanov homology and relations between the Jones polynomial and graph polynomials, we construct a homology theory for embedded graphs from which the chromatic polynomial can be recovered as the Euler characteristic. 
For plane graphs, we show that our chromatic homology can be recovered from the Khovanov homology of an associated link. We apply this connection with Khovanov homology to show that the torsion-free part of our chromatic homology is independent of the choice of planar embedding of a graph.

We extend our construction and categorify the Bollob\'{a}s-Riordan polynomial (a generalisation of the Tutte polynomial to embedded graphs). We prove that both our chromatic homology and the Khovanov homology of an associated link can be recovered from this categorification.

\end{abstract}

\maketitle

\section{Introduction}
\label{sec.intro}

There are numerous connections between graph polynomials and knot invariants in the literature. Perhaps the best known  connection between knot and graph  polynomials is due to  M.~Thistlethwaite. In his seminal paper \cite{Th}, Thistlethwaite proved that the Jones polynomial of an alternating link in $S^3$ can be recovered as  an evaluation of the Tutte polynomial of a plane graph. 
   Thistlethwaite's Theorem was extended by L.~Kauffman in \cite{Ka} where he showed that the Jones polynomial of any link can be obtained as an evaluation of the  signed Tutte polynomial of an edge-signed plane graph (or equivalently the $+-J$ Potts partition function of a plane graph). 
   
More recently, M.~Khovanov constructed a homological generalization of the Jones polynomial. In his influential paper \cite{Kh}, he  constructed a   bigraded homology theory for knots whose graded Euler characteristic is equal to the Jones polynomial.
Khovanov's homology groups are themselves knot invariants and are in fact strictly stronger knot invariants than the Jones polynomial.  Thus Khovanov constructed a homological generalization of the Jones polynomial. 
With Thistlethwaite's Theorem  in mind, it is natural to question whether relations between graph polynomials and the Jones polynomial ``categorify''. That is to ask if one can  construct  a homology theory for graphs with the two properties  that a given graph polynomial arises as the  Euler characteristic of the homology, and  that the Khovanov homology of a link can be recovered from the graph homology of an associated graph, or the graph homology can be recovered from the Khovanov homology of an  associated link.  We will refer to this type of relation as a ``Thistlethwaite-type relation''.  This question on graph and knot homologies motivates the material presented here.   
Before we move on from these motivational considerations,  we  consider additional desirable properties that we would like such a graph homology to have. 

From the point of view of graph theory, Thistlethwaite's connection between the Jones and Tutte polynomials is a little unsatisfactory in  that the relation is between links and plane graphs. We do not want to impose any planarity conditions on our homology theories. We would rather consider homology theories for graphs embedded in  surfaces of any genus. This means that  we would like to construct a homology theory for embedded graphs, such that when the graph is a plane graph, then we obtain the desired relations with Khovanov homology.       
Therefore instead of considering graphs and their polynomials, we consider fatgraphs and their polynomials. A fatgraph is a graph equipped with a cyclic ordering of the incident half-edges at each vertex. Fatgraphs  capture the essential part of an embedded graph.  We note that fatgraphs are also known in the literature as ``ribbon graphs'' and ``maps'', but here    favour the term  ``fatgraph'' which is standard in theoretical physics (see for example \cite{DiFrancesco}).  Rather than working with the Tutte polynomial, when dealing with fatgraphs we instead consider the Bollob\'{a}s-Riordan polynomial \cite{BR1,BR}. This is a  recently defined generalization of the Tutte polynomial to ribbon graphs which captures some of the topology of the fatgraph. 
 Thistethwaite's Theorem relating the Jones polynomial and the Tutte polynomial of a plane graph was recently generalized  by  S.~Chmutov and I.~Pak. In \cite{CP} (which was published in a revised form \cite{CP2}), the 
 Jones polynomial of a (certain type of)  link in a thickened surface was shown to be an evaluation of the Bollob\'{a}s-Riordan polynomial of an associated fatgraph. Furthermore, when the surface is of genus zero, Chmutov and Pak's result specialises to Thistlethwaite's Theorem.
 We note there is currently interest in connections between knots and their polynomials and fatgraphs and their polynomials \cite{CV,Da,HM,Mo,Mo2}.
  This discussion of motivates the categorifications of fatgraph polynomials proposed herein.
 
 \medskip

The paper is structured as follows. After making some preliminary definitions, in Section~\ref{construction} we construct a bigraded chain complex using the set of spanning subfatgraphs  of a fatgraph. We then show that the graded Euler characteristic of the homology of this complex is the chromatic polynomial, thus we have categorified the chromatic polynomial.
Some properties of this homology and connections with other homology theories in the literature   are then given in  Section~\ref{sec:props}. Motivated by the discussion above, in Section~\ref{sec.SBR} we consider categorifications of the  Bollob\'{a}s-Riordan polynomial of a fatgraph. We show how the construction of our chromatic homology can be extended to give a homology theory from which  the Bollob\'{a}s-Riordan polynomial can be  recovered as the graded Euler characteristic. We then prove that both our chromatic homology from Section~\ref{construction}, and  the Khovanov homology of an associated link can be recovered from our fatgraph homology.  This provides our first Thistlethwaite-type relation between graph and knot homology theories.

In Section~\ref{sec:khov} we reconsider our chromatic homology and prove some connections with Khovanov homology. This section contains two main results. One of the main results provides a second Thistlethwaite-type relation which states that our chromatic homology for a plane graph can be recovered from the Khovanov homology of an associated link. The other main result in this section states that the torsion-free part of our chromatic homology is independent of the choice of embedding of a plane graph. The proof of this result utilises the relation with Khovanov homology as well as some recent results on Khovanov homology. 

In the final section we provide a relation between L.~Helme-Guizon and Y.~Rong's categorification of the chromatic polynomial introduced in  \cite{HGR1} and further studied in \cite{CCR,CU,HGR2,HGRP,PPS}, a categorification of the Bollob\'{a}s-Riordan polynomial and Khovanov homology. Specifically we construct a homology theory for the Bollob\'{a}s-Riordan polynomial which comes equipped with two natural homomorphisms: one to    Helme-Guizon and Rong's chromatic homology and the other to Khovanov homology.  
Thus both of these homology theories arise from one homology theory for the Bollob\'{a}s-Riordan polynomial. This addresses the question ``What is the relationship [of Helme-Guizon and Rong's  chromatic homology] with the Khovanov homology for knots?", which was  posed in \cite{HGR1}. 

\medskip

We would like to thank Jo Ellis-Monaghan, Bojan Mohar and Irasema Sarmiento for helpful discussions. M.~L. gratefully acknowledges the support of CONICYT via grant Anillo en Redes.

\section{Some preliminaries}\label{sec:defs}
This section contains some preliminary definitions and results on graphs, fatgraphs, fatgraph polynomials, graded modules. Having set up enough notation, we also provide     a more detailed statement of our results.

\subsection{Fatgraphs}\label{sub.fat}

Let $G=(V,E)$ be an undirected graph, possibly with loops and multiple edges. Each subgraph $(V,W)$, $W\subset E$ of $G$
is called a {\em spanning subgraph}.  Let us denote by $\s(G)$ the set of all spanning subgraphs of $G$.
A graph $F$ is called a {\em fatgraph} if
for each vertex $v\in V$, there is a fixed cyclic order on half-edges adjacent to $v$ (loops are counted twice).
A fatgraph $F$ may be regarded as a  2-dimensional surface with boundary, which will also be denoted by $F$. The
surface is obtained from the fatgraph by fattening the vertices into discs (we will  call these {\em islands})
and connecting them by untwisted fattened edges (which we call {\em bridges}) as prescribed
by the cyclic orders.
The {\em genus}, $g(F)$ of a fatgraph $F$ is defined to be the genus of this surface.
 It will always be clear from the context whether by $F$ we mean the fatgraph or the surface.
We  restrict ourselves  to orientable surfaces.
For a fatgraph $F$ we will usually denote its underlying graph by $G=G(F)$.
Let $V(F)$ be its set of vertices, $E(F)$ its set of edges, and let $v(F)= |V(F)|$,
$e(F)=|E(F)|$, $r(F)= |V|- k(F)$ and $n(F)= |E(F)|- r(F)$.
We denote the number of connected components
of $F$ by $k(F)$, and  the number of connected components
of the boundary of surface $F$ by $p(F)$. The functions $v, e, r, n, k$ will be used for graphs as well. 
Finally, if $F=(V(F),E(F))$ is a fatgraph then each subgraph $F=(V(F),W)$, $W\subset E$ of $F$
is called a {\em spanning fatsubgraph}. We  denote the set of all spanning subgraphs of $G$  by $\s(G)$.

\subsection{Graph polynomials}\label{sub.pol}

Let us recall the definitions of the Tutte and the chromatic polynomials of a graph $G= (V,E)$:
\[
T(G,x,y)= \sum_{H\in \s(G)}(x-1)^{r(G)-r(H)}(y-1)^{n(H)},
\]
\[
M(G,u)= \sum_{H\in \s(G)} (-1)^{e(H)}u^{k(H)}.
\]
The chromatic polynomial $M(G,u)$ is a straightforward evaluation of $T(G,x,y)$.

In \cite{BR1} and  \cite{BR}, Bollob\'{a}s and Riordan defined a fatgraph generalization of the Tutte polynomial.  This  three variable polynomial is defined by the state sum
\begin{equation}\label{eq:BRpoly}
R(F,x,y,z)  = \sum_{H\in \s(F)}x^{r(F)-r(H)}y^{n(H)}z^{k(H)-p(H)+n(H)}.
\end{equation}
The exponent of $z$ is equal to twice the genus of the fatgraph $H$ and we may therefore write
\[
R(F,x,y,z) = \sum_{H\in \s(F)}x^{r(F)-r(H)}y^{n(H)}z^{2g(H)}.
\]


 If $F$ is a fatgraph and $G$ its underlying graph,  one can express
the chromatic polynomial of $G$ in terms of geometric information from $F$.
The next lemma expresses the chromatic polynomial in the evaluation we  use.
\begin{lemma}\label{lm:cprw2}
Let $F$ be a fatgraph and $G$ be its underlying graph. Then
\begin{multline}\label{eq:cprw2}
(q+q^{-1})^{e(F)}M(G, (q+q^{-1})^2) = \\
(-1)^{e(F)}(q+q^{-1})^{v(F)}\sum_{H\in \s(F)}(q+q^{-1})^{(p(H)+2g(H))}(-q)^{e(F)-e(H)}(1+q^{-2})^{e(F)-e(H)}.
\end{multline}
\end{lemma}
\begin{proof}
Using the identity $2g(H)= k(H)-p(H)+n(H)$  and the definitions of $n(H),r(H)$ above, we have
$2k(H)= p(H)+2g(H)-e(H)+v(H)$. Hence
\[
M(G,u)= u^{1/2v(F)}\sum_{H\in \s(F)}u^{1/2p(H)+g(H)}[(-1)u^{1/2}]^{-e(H)}.
\]
Substituting $u^{1/2}= (q+q^{-1})$ we get
\[
M(G,(q+q^{-1})^2)= (q+q^{-1})^{v(F)}\sum_{H\in \s(F)}(q+q^{-1})^{p(H)+2g(H)}(-1)^{-e(H)}(q+q^{-1})^{-e(H)}.
\]
This easily implies Equation~\ref{eq:cprw2}.
\end{proof}

Our main  object of  study is the above evaluation and scaling of the chromatic polynomial. We set
\begin{multline*} 
  Z(F,q)    = (q+q^{-1})^{e(F)} M(G, (q+q^{-1})^2)  \\
   =   (-1)^{e(F)}(q+q^{-1})^{v(F)} \sum_{H\in \s(F)} (q+q^{-1})^{p(H)+2g(H)}
(-q)^{e(F)-e(H)}(1+q^{-2})^{e(F)-e(H)} .
  \end{multline*}

\subsection{Graded modules}\label{sub.results}

Let $M=\bigoplus_{i \in \mathbb{Z} }M_i$ be a graded $\mathbb{Z}$-module. The \emph{graded dimension} of $M$ is defined by
\[
\text{qdim}(M) :=  \sum_{i \in \mathbb{Z}} q^i \text{rk} \left( M_i \right)
= \sum_{i \in \mathbb{Z}} q^i \text{dim}_{\mathbb{Q}} \left( M_i \otimes_{\mathbb{Z}} \mathbb{Q} \right).
\]
If $H = \left( H^i \right)_{i \in \mathbb{Z} }$ is the homology of some chain complex of graded $\mathbb{Z}$-modules, the {\em Poincar\'{e} polynomial} is the two-variable Laurent polynomial
\[
P(H) = \sum_{i \in \mathbb{Z}} t^i \qdim\left( H^i \right) \in \mathbb{Z}[q,q^{-1},t,t^{-1}].
\]
The Poincar\'{e} polynomial encodes all of the torsion-free information of the homology groups.
The {\em Euler characteristic} is defined to be the evaluation $ \chi (H) = P(H)(t=-1)$. It  generalizes  the usual
Euler characteristic of graphs and surfaces.

\medskip

We construct homology theories for fatgraphs which have the property that a given graph polynomial can be recovered as its Euler characteristic. 
For the convenience of the reader we summarise the main results of this paper in the following theorem.
\begin{theorem}\label{thm.one}

Let $F$ be a fatgraph and $G$ be its underlying graph. Let $\C (F)$ be its chain complex as constructed in Section \ref{construction}. Then the following hold:
\begin{enumerate}
\item The Euler characteristic of the homology ($H^i(\C (F)) := (\ker  d^i)/({\rm im} \, d^{i-1}) $) is equal to
the chromatic polynomial $Z(F,q)$.

\item The homology groups $H$ are strictly stronger graph invariants than the chromatic polynomial.

\item The Poincar\'{e} polynomial is invariant on different planar embeddings of a planar graph $G$. However, the homology is dependent upon the genus of the embedding of a graph.

\item The chromatic homology of a plane graph can be recovered from the  Khovanov homology of an associated link.

\item This homology theory may be extended to a categorification
of the Bollob\'{a}s-Riordan polynomial of the signed fat graphs, from which
the Khovanov homology of an associated link may be recovered. Our chromatic homology can also be recovered from this homology.
\end{enumerate}
\end{theorem}
This theorem will follow from Theorems~\ref{th:homology}, which contains statement 1; Proposition~\ref{basics}, which gives some properties of the homology;
 Theorem~\ref{th:invariance} which contains statement 3, Theorem~\ref{th:univcoeff}   contains statement 4 and
section \ref{sec.SBR} which contains the construction for statement 5 above.

\section{Construction of the homology} \label{construction}

Let  $F$ be a fatgraph and $G=(V,E)$ be its underlying graph.
We call a spanning fatsubgraph of a fatgraph a {\em state}. The chromatic polynomial $Z(F,q)$ is expressed in (\ref{eq:cprw2}) as a sum over all states.
Each state of a fatgraph $F$ is obtained by the removal of  a set of bridges of $F$. For example:
\[
\begin{array}{ccc}
\includegraphics[width=4cm]{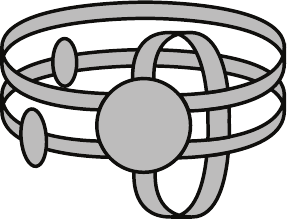}
 & \hspace{1.5cm} & 
 \includegraphics[width=4cm]{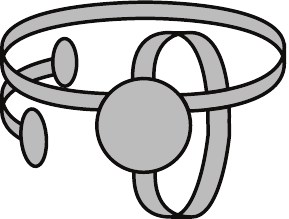}
\\
\text{{\em A fatgraph.}} & &\text{{\em A state.}}
\end{array}.
\]
We call the total number of edges of $F$ minus the number of bridges in a state $H$ the {\em height} of $H$, denoted by $h(H)$, so $h(H)=e(F)-e(H)$.

Let us consider the following example of the state sum $Z(F,q)$:
\[\includegraphics[width=15cm]{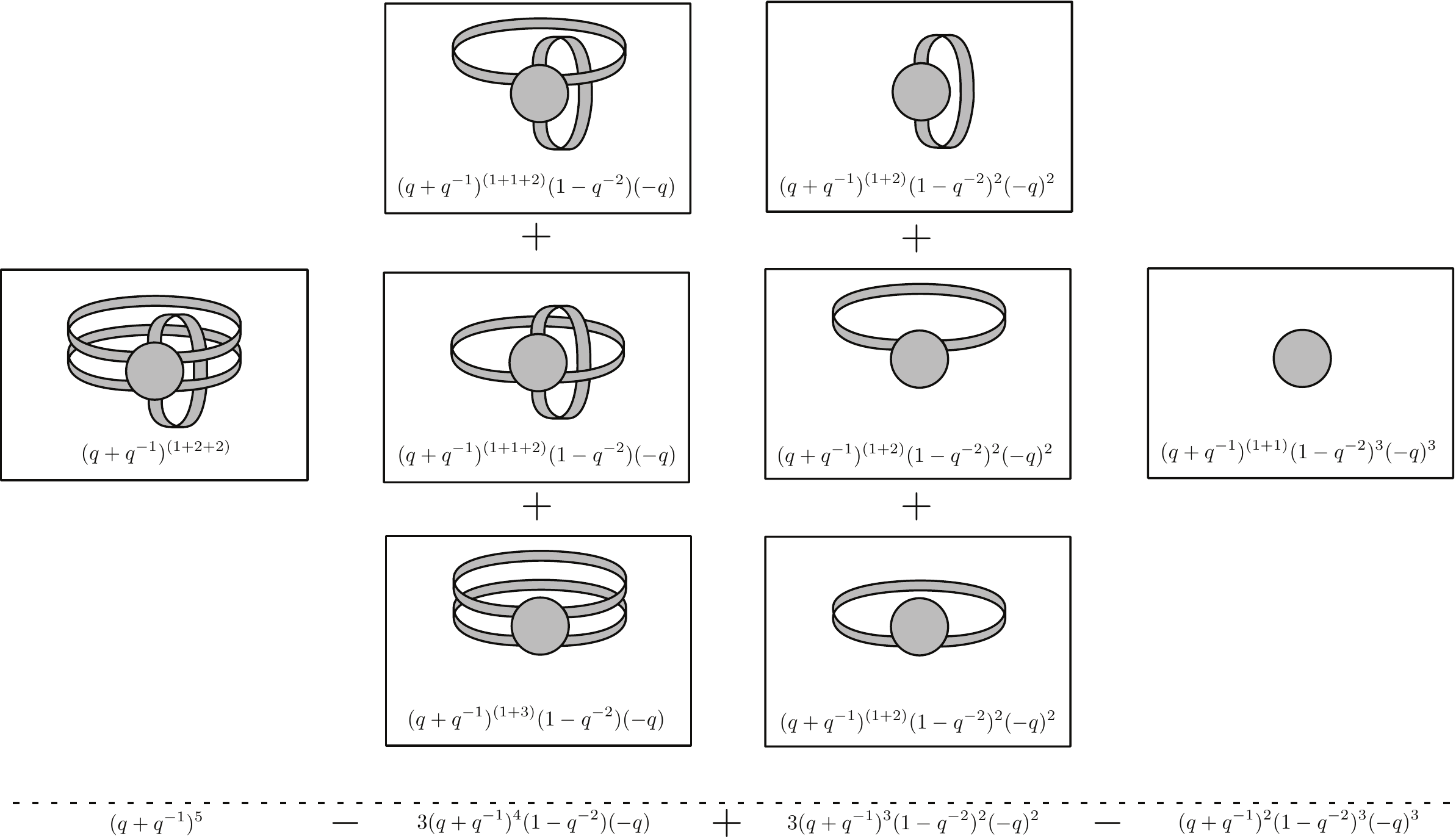}
, \]
which, of course, is equal to zero.
We will use this example to illustrate how the state sum for $Z(F,q)$ gives rise to a chain complex.
The approach taken here is similar to Bar-Natan's exposition of Khovanov homology in \cite{BN1}.

Notice that each state of the fatgraph in Equation~\ref{eq:cprw2} gives rise to (up to sign) a Laurent polynomial of the form
\begin{equation} \label{eq:poly}
(q+q^{-1})^{(v+p+2g)} ( 1+q^{-2})^h q^h,
\end{equation}
where $v=v(H)$, $p=p(H)$, etc..
We want to replace each such polynomial term of the state sum $Z(F,q)$ with a graded module whose graded dimension is equal to this polynomial.

To do this we define $V$ to be the free, graded $\Z$-module with two basis elements $v_-$ and $v_+$ in graded degrees $-1$ and $+1$ respectively, and $R$ to be the free, graded $\Z$-module with basis elements $x_{-2}$ and $x_0$ in graded degrees $-2$ and $0$. Notice that $\qdim (V) = q+ q^{-1}$ and $\qdim (R) = 1+ q^{-2}$.

The {\em degree shift} operation $\{ \cdot \}$ on graded modules is defined by setting \[ M\{ l\}_m := M_{m-l}.\] Clearly $\qdim (M\{ l\}) = q^l \cdot  \qdim (M)$.  We  note that $R = V\{-1\}$ and therefore every occurrence of $R$ in this paper could be replaced with $V\{-1\}$. We retain the  use of $R$ for clarity, however the reader should bear in mind that the two modules only differ by a grading shift. We will see that, in some sense,  $R$ plays the role of coefficients in the homology theory.

Upon observing that for two graded modules $M$ and $N$,
\[ \qdim (M \otimes N) =  \qdim (M ) \cdot \qdim (N) \] and \[
\qdim (M \oplus N) =  \qdim (M ) + \qdim (N),\]
 it is easily seen that the modules
\begin{equation} \label{eq:state}
 V^{\otimes (v+p+2g)} \otimes R^{\otimes h} \{h\}
\end{equation}
have  graded dimensions equal to the Laurent polynomial (\ref{eq:poly}).

In order to simplify the text, we abuse notation and identify the state of a fatgraph with its assigned   polynomial term (\ref{eq:poly}) and its assigned graded module (\ref{eq:state}).

Next, as in the calculation for $Z(F,q)$ in the example above, we arrange the states into $n=e(F)$ columns indexed by the height of the corresponding state of the fatgraph, so the $i$-th column contains all modules which come from states of height $i$.
We then define the $i$-th chain module $\tC^i(F)$ to be the direct sum of all of the modules corresponding to states of height $i$.
For example (ignoring the maps in the figure for the time being), the above example of $Z(F,q)$ becomes:
\[\includegraphics[width=15cm]{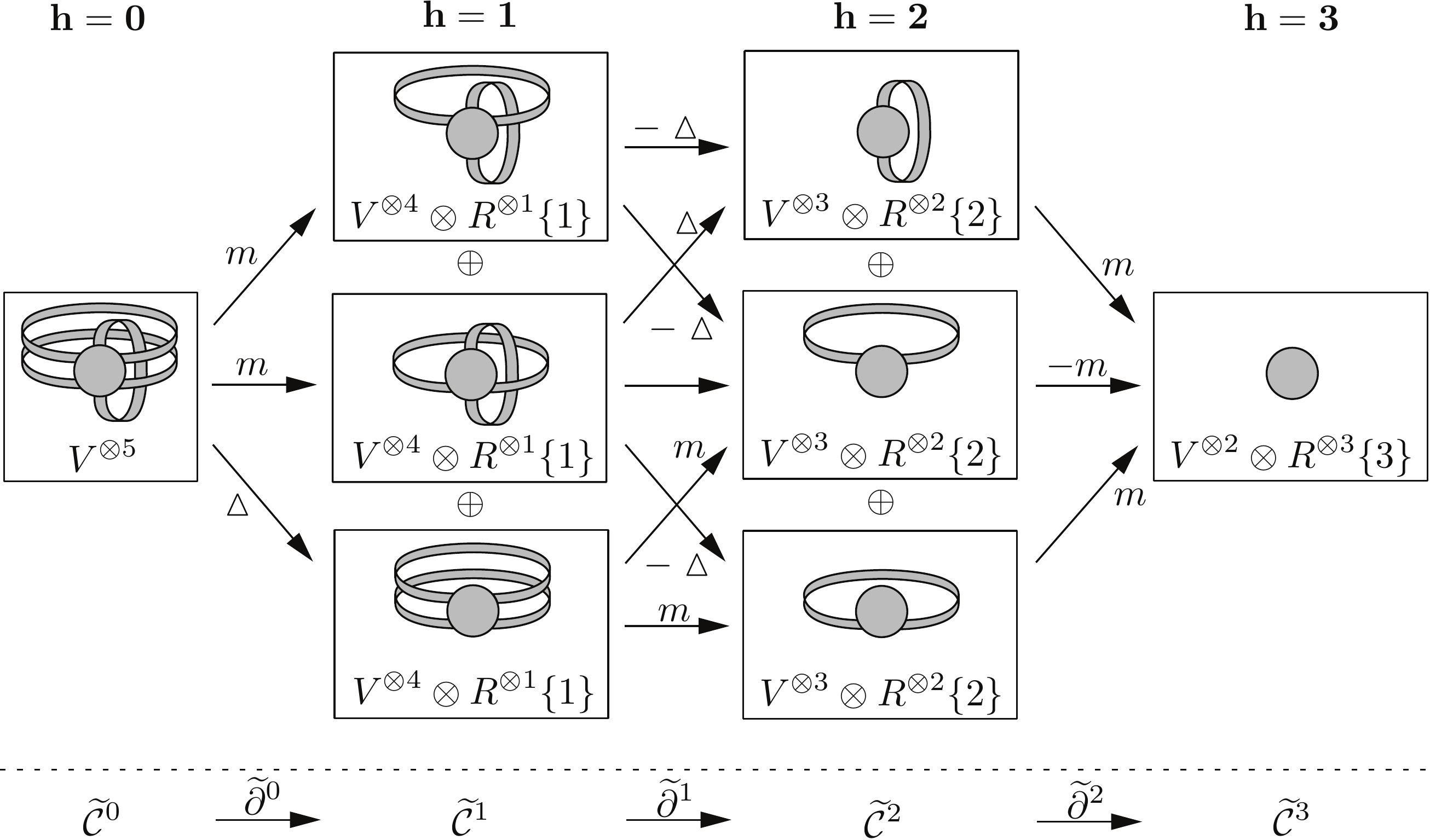}
.\]
Note that although the tensor powers of $V$ are monotone in this example, this is not true in general.

Our next task is to define chain maps.
We begin with the observation that removing a bridge from a state of height $h$ determines a state of height $h+1$. The state  obtained depends on which bridge is removed. Notice that the state determined by the removal of a bridge will have one more or one less boundary component than the original state,  and the genus of the state will either be unchanged or will decrease by one.
Whenever we can move from a state of height $h$ to a state of height $h+1$ by the removal of an edge we will define a {\em per-edge map}
\[
\delta : V^{\otimes (v+p+2g)} \otimes R^{\otimes h} \{h\}
\longrightarrow
V^{\otimes (v+p^{\prime}+2g^{\prime})} \otimes R^{\otimes (h+1)} \{h+1\},
\]
where $p^{\prime} = p \pm 1$ and $g^{\prime} = g$ or $g-1$.
In order to  define $\delta$, we first describe three maps which correspond to the  tensor factors $V^{\otimes p}$, $V^{\otimes 2g}$ and $R^{\otimes h}$.

\medskip

\noindent \underline{\textbf{$V^{\otimes p}$:}}
In terms of the boundary components of the states, one of two things can happen when we move from one state to another by the deletion of a bridge: either  two boundary components will be merged into one component or one component will be split into two components. For example
\[
\includegraphics[height=15mm]{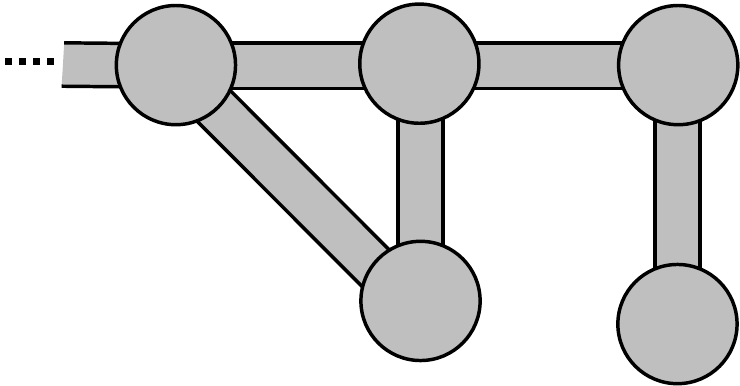}
\hspace{5mm}
\raisebox{6mm}{\includegraphics[width=7mm]{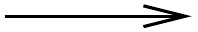}}
\hspace{5mm}
\includegraphics[height=15mm]{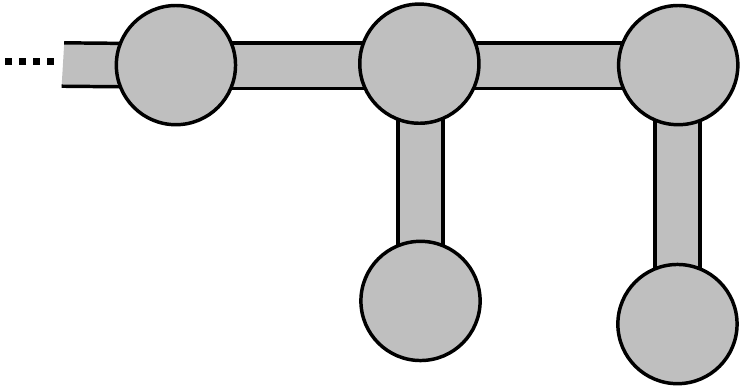}
\hspace{5mm}
\raisebox{6mm}{\text{or}}
\hspace{5mm}
\includegraphics[height=15mm]{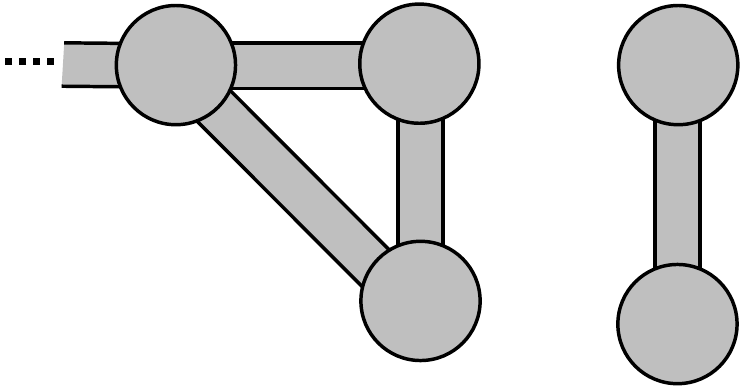}
\raisebox{6mm}{.}
\]
By (\ref{eq:state}), a copy of the module $V$ is assigned to each boundary component of a state. We fix a correspondence between the boundary components and the modules $V$ once and for all (the homology will be independent of this choice). When two components merge we need a multiplication map \[ \widetilde{m}: V^{\otimes p} \rightarrow V^{\otimes (p-1)}.\] We define the map to be  the identity on all copies of $V$ which are unchanged by the addition of an edge, and to act on the two merging components by the multiplication defined on basis elements by
\[
\widetilde{m}^{\prime}:
\left\{
\begin{array}{lll}
v_- \otimes v_- \mapsto 0 & & v_+ \otimes v_- \mapsto v_- \\
v_+ \otimes v_+ \mapsto v_+ & & v_- \otimes v_+ \mapsto v_-
\end{array}\right.
.
\]
In the case where one component is split into two, we define a coproduct \[\widetilde{\vartriangle}: V^{\otimes p}  \rightarrow V^{\otimes (p+1)}\] which is the identity on all factors of $V^{\otimes p}$ except on  the component being split where it acts on basis elements by
\[
\widetilde{\vartriangle}^{\prime}:
\left\{
\begin{array}{l}
v_+ \mapsto  v_+ \otimes v_-  +  v_- \otimes v_+  \\
 v_- \mapsto v_- \otimes v_-
\end{array}\right.
.
\]
Notice that the two maps  $\widetilde{m}^{\prime}$ and $\widetilde{\vartriangle}^{\prime}$ are the maps used in the definition of Khovanov homology (\cite{Kh, BN1}). This observation will prove to be important later when we find Thistlethwaite-type relations between graph and knot homologies.
\medskip

\noindent \underline{\textbf{$V^{\otimes 2g}$:}}
We identify $V^{\otimes 0}$ with $\mathbb{Z}$.
There are two cases. If $g=g^{\prime}$ then we set $\widehat{m}_g: V^{\otimes 2g} \rightarrow V^{\otimes 2g}$ equal to $1^{\otimes 2g}$.
If $g^{\prime}=g-1$ then we define
$\widehat{m}_g: V^{\otimes 2g} \rightarrow V^{\otimes 2(g-1)}$ on its basis elements by
\[
\widehat{m}_g :
v_{i_1} \otimes \cdots \otimes v_{i_{2g}}
\mapsto
\sum v_{j_1} \otimes \cdots \otimes v_{j_{2g-2}},
\]
where the sum is over all basis elements of $V^{\otimes 2(g-1)}$ whose graded dimension is equal to the graded dimension of the basis element $v_{i_1} \otimes \cdots \otimes v_{i_{2g}}$ of $V^{2g}$.  $\widehat{m}_g$ is a map of graded degree $0$.

\medskip

\noindent \underline{\textbf{$R^{\otimes h}$:}}
Setting $R^0 = \mathbb{Z}$, we define $\overline{\vartriangle}_h : R^{\otimes h} \rightarrow R^{\otimes (h+1)} $ by
\[
\overline{\vartriangle}_h :
\left\{
\begin{array}{lll}
1 \mapsto x_0 , & & \text{when } h=0 \\
y \otimes x_i \mapsto  y \otimes \sum_{ k+l=i} x_k \otimes x_l , &  & \text{otherwise}
\end{array}\right. .
\]
Again this is a graded degree $0$ map.

\medskip
We take the tensor of these maps and define maps 
\[ V^{\otimes (v+p+2g)} \otimes R^{\otimes h} \{h\} \longrightarrow V^{\otimes (v+p^{\prime}+2g^{\prime})} \otimes R^{\otimes (h+1)} \{h+1\}\]
 by
\[ \begin{array}{ccc}
m = 1^{\otimes v}  \otimes  \widetilde{m} \otimes  \widehat{m}_g \otimes \overline{\vartriangle}_h & \text{and} & \vartriangle = 1^{\otimes v}  \otimes  \widetilde{\vartriangle} \otimes  \widehat{m}_g \otimes \overline{\vartriangle}_h
\end{array}.  \]
We need to assign a sign $+1$ or $-1$ to each of the maps $m$ and $\vartriangle$ to obtain the per-edge maps $\delta$. This is done as in \cite{Kh, BN1} by realizing the states as  vertices of an $n$-dimensional cube and the per-edge maps as its edges, where $n=|E|$. To do this  we label the bridges of $F$ with $1, \ldots , n$. The homology is independent of the choice of  labeling. A  proof of this follows the proof of \cite{HGR1} Theorem~2.12 and is therefore excluded. Each state of $F$ can be represented by the vertex of a $n$-dimensional cube $( \alpha_1 , \ldots , \alpha_n )$, by setting $\alpha_i = 0$ if the bridge labelled $i$ is  in the state, and $\alpha_i = 1$ if it is not.
A per-edge map is a map from a state labelled
$( \alpha_1 , \ldots , \alpha_{j-1}, 0 ,\alpha_{j+1}, \ldots  , \alpha_n )$
 to one labeled
$( \alpha_1 , \ldots , \alpha_{j-1}, 1 ,\alpha_{j+1}, \ldots  , \alpha_n )$.
The sign $(-1)^{\sum_{i<j}\alpha_i}$ is then assigned to each of the maps $m$ and $\vartriangle$ in the complex. This defines the per-edge maps $\delta$.

Finally the differential $\widetilde{\partial}^h : \tC^h \rightarrow \tC^{h+1}$ is obtained as the sum of all of the per-edge maps between the tensor factors of $\tC^h$ and $\tC^{h+1}$.

For a fatgraph $F$, we let $\tC (F)$ denote the complex $( \tC^h ,\widetilde{\partial^h})$ constructed as above.
\begin{lemma}
$\tC$ is a chain complex (\emph{i.e.} $\widetilde{\partial} \circ \widetilde{\partial} = 0$) and the differentials are of graded degree $0$.
\end{lemma}
\begin{proof}
The per-edge maps are easily seen to be  (co)associative and (co)commutative. The first statement then follows as the per-edge maps around each state anti-commute. The second statement follows since $ \widetilde{m}$ and $\widetilde{\vartriangle}$ are of degree $-1$ and are therefore  of degree zero once the target is shifted by $\{1\}$,  and $\widehat{m}_g $ and  $\overline{\vartriangle}_r$ are of degree $0$. 
\end{proof}

Finally, in order to deal with the factor $(-1)^{e(F)}$ of $Z(F,q)$, we define the
 {\em height shift} operation $[ \cdot ]$ on chain complexes by
$ \left( \C^i , \partial^i  \right) [ s ] :=  \left( \C^{i-s} , \partial^{i-s}  \right) $.
We can then normalize the chain complex by  $[-e(F) ]$ and define
\[
\C =(\C^h , \partial^h ) = \tC (F)[ -e(F) ] =( \tC^h ,\widetilde{\partial^h} )[ -e(F)]
. \]

Recall that the \emph{homology} of a chain complex $\C =(\C^h , \partial^h)$ is the sequence $H(\C) = \left( H^i (\C) \right)_{i \in \mathbb{Z}}$, where $H^i (\C)= \ker ( \partial^i )   /  \textrm{im}( \partial^{i-1} )$.

\begin{theorem} \label{th:homology}
Let $F$ be a fatgraph, $G$ its associated graph and $\C (F)$ its chain complex. Then the Euler characteristic of the  homology $H(\C(F))$ is equal to the scaled chromatic polynomial $Z(F,q)$. Moreover the homology is an invariant of fatgraphs and is a strictly stronger invariant than the chromatic polynomial $M(G,q)$.
\end{theorem}
\begin{proof}
The first statement follows essentially by construction.
It is well known that when the differentials are of graded degree $0$, $\qdim \left(H^i(\C)\right)$ is equal to $\qdim \left(\C^i\right)$. In turn this is equal to the sum of the graded dimensions of the states of height $i$. The result then follows. (This is the graded extension of a classical result in topology, see eg \cite{Ha} page~146.) The shift $[-e(F) ]$ ensures that the alternating sum of $\qdim \left(H^i(\C)\right)$ has the correct sign.

The second statement follows by a calculation (see Proposition~\ref{basics}). 
\end{proof}

\section{Properties of the homology}\label{sec:props}
Having having constructed our categorification of the chromatic polynomial, we show that the homology groups satisfy various desirable properties.

\begin{proposition} \label{basics}
The following hold for the homology $H\left(\C (F) \right)$:
\begin{enumerate}
\item The homology groups are strictly stronger than the chromatic polynomial.

\item The homology will differentiate between graphs which differ only by multiple edges or loops.

\item The homology is not a Tutte-Grothendieck  invariant.

\item If $F$ is the disjoint union of fatgraphs $F_1$ and $F_2$, then
\begin{multline*}
H^i\left( \C (F) \right) = \left( \bigoplus_{p+q=i}
 H^p\left( \C (F_1) \right) \otimes H^q\left( \C (F_2) \right)   \right)\\
\bigoplus
\left( \bigoplus_{p+q=i-1} \rm{Tor}_1 \left(   H^p\left( \C (F_1) \right) ,  H^q\left( \C (F_2) \right) \right) \right).
\end{multline*}

\item Let $e$ be a bridge of a fatgraph $F$ and let $F-e$ denote $F$ with the bridge $e$ deleted and $F/e$ denote $F$ contracted along the edge $e$,  then there exists a deletion-contraction exact sequence
\[
\begin{diagram}
H^{\ast}\left(  \C (F - e)  \right)\otimes R & & \rTo  && H^{\ast}\left(  \C (F) \right)  \\
& \luTo  & &\ldTo &\\
& & H^{\ast}\left(  \C (F/e) \right) \otimes V&&
\end{diagram}.
\]

\item If $F_1$ is a subfatgraph of $F_2$ the inclusion map  induces a homomorphism $ H(i): H(F_1) \rightarrow H(F_2)   $ in homology.
\end{enumerate}
\end{proposition}

\begin{proof}

\emph{1.} This follows by a calculation. For example it is easy to check that the fatgraphs
\includegraphics[width=7mm]{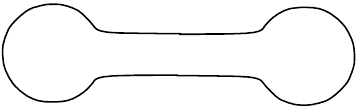}
and \includegraphics[width=7mm]{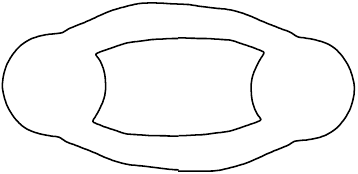}
have the same chromatic polynomials but different chromatic homology.

\emph{2.} This follows since the homology in  the highest degree is non-zero (as $\bar{\vartriangle}_{n-1}$ is not surjective).

\emph{3.} Again this follows by a calculation. For example the fatgraphs
\includegraphics[width=7mm]{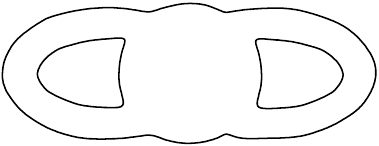}
 and  
 \includegraphics[width=7mm]{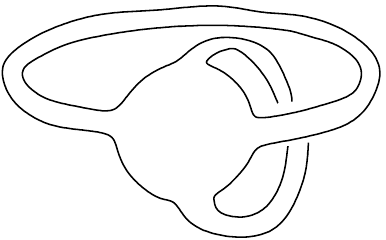}
  have different homology. Note that this example also shows that the Poincar\'{e} polynomial is not a Tutte-Grothendieck invariant.

\emph{4.} As in \cite{HGR1}, the chain complex of $F$ can be written $ \C (F) = \C (F_1) \otimes \C (F_2)$. Then, since the chain complexes are free, the result follows by  an application of the K\"{u}nneth formula (see eg. \cite{Br}).

\emph{5.} Since the homology is independent of the labelling of the bridges,
we may assume that $e$ is the bridge which is labelled last so that $(\alpha_1 , \ldots , \alpha_{n-1}, 0)$ is the vertex label of a state  containing the bridge $e$, while $(\alpha_1 , \ldots , \alpha_{n-1}, 1)$ is the label of the corresponding state with the bridge $e$ deleted.
Let $\alpha = (\alpha_1 , \ldots , \alpha_{n-1})$ be a state of $F-e$ of height $h$. Then $\alpha^{\prime} = (\alpha_1 , \ldots , \alpha_{n-1},1)$ gives  a state of $F$.
Notice that the fatgraphs of the states $\alpha$ and  $\alpha^{\prime}$ are identical but  the state $\alpha^{\prime}$ is of height $h+1$. Therefore
  $V^{\otimes (v+p+2g)} \otimes R^{\otimes h} \{ h \}$ is the module assigned to the state $\alpha$. The map
\[
\eta_{\alpha} : V^{\otimes (v+p+2g)} \otimes R^{\otimes h} \{ h \} \otimes R
\longrightarrow  V^{\otimes (v+p+2g)} \otimes R^{\otimes (h+1)} \{ h +1 \},
\]
which takes basis elements to themselves, then induces a map \[ \eta : \C (F-e) \otimes R \rightarrow \C (F).\]

We also need to define a map $\nu : \C (F) \rightarrow \C (F/e) \otimes V$.
To do this, suppose that $\alpha^{\prime} = (\alpha_1 , \ldots , \alpha_{n})$  is assigned to a  state of $F$. Then
$\alpha^{\prime \prime} = (\alpha_1 , \ldots , \alpha_{n-1})$ is a state of $F/e$. If $\alpha_{n} =0$, then the fatgraph in the state $\alpha^{\prime \prime}$ is equal to the fatgraph of $\alpha^{\prime}$ contracted along the bridge $e$. Notice that the height of the two states $\alpha^{\prime }$ and  $\alpha^{\prime \prime}$ are  equal, however the fatgraphs $F$ has one more island than $F/e$. We then define a map
\[\nu_{\alpha^{\prime}} : V^{\otimes (v+p+2g)} \otimes R^{h} \{h \}   \longrightarrow V^{\otimes ((v-1)+p+2g)} \otimes R^{h} \{h \} \otimes V \] by
\begin{multline*}
\nu_{\alpha^{\prime}} :
 v_1 \otimes \cdots v_{v+p+2g} \otimes x_1 \otimes \cdots \otimes x_h \\
\longmapsto
\left\{
\begin{array}{lll}
  v_2 \otimes \cdots  \otimes x_h \otimes v_1   & & \text{if } \alpha^{\prime} = ( \alpha_1 , \ldots , \alpha_{n-1},0 ) \\
0  &  & \text{otherwise}
\end{array}\right. ,
\end{multline*}
where $v_1 \otimes \cdots v_{v+p+2g} \otimes x_1 \otimes x_h$ denotes a basis element of $V^{\otimes (v+p+2g)} \otimes R^{h} \{h \}$.
These  induce $\nu : \C (F) \rightarrow \C (F/e) \otimes V$.

The next step is to show that $\eta$ and $\nu$ are chain maps, that is  $\partial \circ \eta = \eta \circ \partial$ and   $\partial \circ \nu = \nu \circ \partial$.
It is enough to show this for the per-edge maps, \emph{i.e.} that $\eta$ and $\nu$ commute with $\delta$.
Since each $\eta_{\alpha}$ maps a state to a state consisting of  the same fatgraph, the corresponding per-edge maps in the two complexes will both be of the form of $\pm m$ or $\pm \vartriangle$, and since we chose the edge $e$ to be the last element $\alpha_n$ in the vertex labeling, the maps will occur with the same sign in both complexes. Therefore it is enough to show $m \circ \eta = \eta \circ (m \otimes 1)$ and  $\vartriangle \circ \eta = \eta \circ (\vartriangle \otimes 1)$, which is easily verified.

Similarly, since $\nu$ maps a fatgraph to $0$ or the fatgraph contracted along the edge $e$, it is enough to show
$m \circ \nu = \nu \circ (m \otimes 1)$ and  $\vartriangle \circ \nu = \nu \circ (\vartriangle \otimes 1)$. Again this is easily verified.

It is not hard to see that the chain maps $\eta$ and $\nu$ form a short exact sequence of chain complexes:
\[
0 \rightarrow \C^{i-1} (F-e) \otimes R \rightarrow \C^{i} (F) \rightarrow  \C^{i} (F/e) \otimes V \rightarrow 0 .
\]
This induces the long exact sequence in homology, completing  the proof.

\emph{6.}
There is a natural inclusion $F_1 \hookrightarrow F_2$. This extends to an inclusion map between states: for each state $H$ of $F_1$ there is a state $H_2$ of $F_2$ whose fatgraph consists of $H$ and $k$ independent islands,
 where $k=v(F_2) - v(F_1)$.
This inclusion between states of fatgraphs induces  an inclusion
\[ V^{\otimes v(H)+p(H)+2g(H)} \otimes R^{\otimes h(H) } \{ h (H)\} \hookrightarrow
  V^{\otimes (v(H_2))+(p(H_2))+2g(H_2)} \otimes R^{\otimes h(H_2) } \{ h(H_2) \}  \]
between modules and therefore a map $\C(F_1) \hookrightarrow \C(F_2)$ of chain complexes.
This is clearly a chain map  and therefore induces a homomorphism in homology.
(Note that this property also holds for the Tutte homology constructed in  \cite{JHR}.) 
\end{proof}

\begin{remark}
The homology theory is genuinely different from the categorification of the chromatic polynomial for abstract graphs constructed in \cite{HGR1}. This can be seen immediately from properties \emph{2.} and \emph{3.} of Proposition~\ref{basics}. In fact  Helme-Guizon and Rong's  homology groups are trivial on graphs which contain loops and can not distinguish graphs which differ only by multiple edges (see \cite{HGR1}).
\end{remark}

\begin{remark}
As with \cite{HGR2}, our construction works over a variety of algebras. However, Khovanov's homology of links requires a unique choice of algebra  for isotopy invariance, and since our proof of the invariance of the Poincar\'{e} polynomial for  planar fatgraph homology (Theorem~\ref{th:invariance}) is dependent upon the invariance of Khovanov homology, we do not know if the corresponding result holds over  different algebras.
\end{remark}

As mentioned in the introduction, our homology theory is motivated by Khovanov's categorification of the Jones polynomial and Thistlethwaite's Theorem connecting the Jones and Tutte polynomials. We were motivated by the problem of lifting Thistlethwaite's Theorem to the homological level. 
For  planar fatgraphs, such a relation will be given later in Theorem \ref{th:univcoeff}.
This theorem states that for fatgraphs
of genus zero, our homology is obtained by ``adding coefficients'' to the Khovanov's homology
of the corresponding alternating link diagram. 
Before we discuss this connection with Khovanov homology, we will consider homology theories for the Bollob\'{a}s-Riordan polynomial of a fatgraph.

\section{The Bollob\'{a}s-Riordan Polynomial}
\label{sec.SBR}

A promised in the introduction, we provide  categorifications of the Bollob\'{a}s-Riordan polynomial of a fatgraph. In Subsection~\ref{ss:sbr} we consider the categorification of a one-variable specialization of the Bollob\'{a}s-Riordan polynomial. We then state a proposition which relates this Bollob\'{a}s-Riordan homology with the chromatic homology considered above. This relationship between the two fatgraph homologies will prove to be important in Section~\ref{sec:khov} where a connection between our chromatic homology and Khovanov homology is given. 

In Subsection~\ref{sec.MBRP} we go on to extend our categorification of the  one-variable specialization of the Bollob\'{a}s-Riordan polynomial to a categorification of the full three-variable homology by considering trigraded modules. We conclude the section with a proof that Khovanov homology and our chromatic homology can be recovered from the categorification of the Bollob\'{a}s-Riordan polynomial. This provides one of our homological analogues of Thistlethwaite's Theorem.

\subsection{Chromatic homology as an extension of Bollob\'{a}s-Riordan homology}\label{ss:sbr}
Here we categorify a one-variable specialization of the Bollob\'{a}s-Riordan polynomial. To do this we first write a suitable evaluation of $R(F,x,y,z)$ as a state sum.
\begin{lemma}\label{lm:br5}
Let $F$ be a fatgraph and $G$ be its underlying graph. Then
\begin{multline}\label{eq:br5}
\left. \left( x^{k(F)}y^{v(F)}[-y^{-1}(xy)^{1/2}]^{e(F)}R(F,x,y,z)\right)\right|_{
x=  -q(q+q^{-1}),\, y=  -q^{-1}(q+q^{-1}),\, z=1} \\
= (-1)^{e(F)}\sum_{H\in \s(F)}(q+q^{-1})^{v(H)+p(H)+ 2g(H)}(-q)^{e(F)-e(H)}.
\end{multline}
\end{lemma}
\begin{proof}
Using $r(F)= v(F)- k(F)$, $n(F)= e(F)- r(F)$ and $g(F)= 1/2(k(F)-p(F)+n(F))$, we have
\begin{equation*}
\begin{split}
R(F,x,y,z) & = x^{-k(F)}y^{-v(F)}\sum_{H\in \s(F)}(xy)^{k(H)}y^{e(H)}z^{2g(H)} \\
& = x^{-k(F)}y^{-v(F)}\sum_{H\in \s(F)}(xy)^{1/2(v(H)+p(H)+ 2g(H))}[y(xy)^{-1/2}]^{e(H)}z^{2g(H)}
\\
& = x^{-k(F)}y^{-v(F)} [y^{-1} (xy)^{1/2}]^{-e(F)}
\\ & \hspace{18mm} \sum_{H\in \s(F)}(xy)^{1/2(v(H)+p(H)+ 2g(H))}[y^{-1}(xy)^{1/2}]^{e(F)-e(H)}z^{2g(H)}
.
\end{split}
\end{equation*}
The lemma follows upon substituting $x=  -q(q+q^{-1})$, $y=  -q^{-1}(q+q^{-1})$ and $z=1$. 
\end{proof}

We define the {\em restricted Bollob\'{a}s-Riordan polynomial} by
$$
\widehat{R}(F,q)= \sum_{H\in \s(F)}(q+q^{-1})^{v(H)+p(H)+ 2g(H)}(-q)^{e(F)-e(H)}.
$$
By assigning the modules \[ V^{\otimes (v(H)+p(H)+ 2g(H))}  \{ h \} \] to a state which contributes the summand \[ (q+q^{-1})^{v(H)+p(H)+ 2g(H)}(-q)^{h}\] to $\widehat{R}$, we may categorify $\widehat{R}$  using a construction similar to that used  in Section~\ref{construction} to categorify the chromatic polynomial.
The differentials are defined through a restriction of the per-edge maps used in the above complex for the chromatic polynomial $Z(F,q)$, so that
\[m_{\widehat{R}}=  1^{\otimes v(H)} \otimes \widetilde{m} \otimes  \widehat{m}_g \] and
\[ \vartriangle_{\widehat{R}} = 1^{\otimes v(H)} \otimes  \widetilde{\vartriangle} \otimes  \widehat{m}_g.\]
One then obtains a chain complex $\hC(F)$ which we call the {\em restricted Bollob\'{a}s-Riordan chain complex}.
We denote the homology of this complex by $\hH(\hC(F))$. Just as before we have:
\begin{theorem} \label{th:BRhomology}
Let $F$ be a fatgraph and $\hC (F)$ be its restricted
Bollob\'{a}s-Riordan chain complex. Then the Euler characteristic of the  homology $\hH(\hC(F))$ is equal to the restricted Bollob\'{a}s-Riordan polynomial $\widehat{R}(F)$.
\end{theorem}

The following universal coefficient type theorem relates the categorification of the chromatic polynomial from Section~\ref{construction} to this categorification of the restricted Bollob\'{a}s-Riordan polynomial $\widehat{R}(F)$.
Although the proposition can be proved directly, it will also follow from  Proposition~\ref{th:UCT} which is stated and proved in  the following subsection. Consequently we prefer to delay the proof Proposition~\ref{th:UCT} until we can prove it as an application of the stronger theorem.
\begin{proposition} \label{th:BRuniv}
Let $F$ be a  fatgraph. Then
\[
\widetilde{H}^i_j\left( \widetilde{\C} (F) \right)
 = \bigoplus_{p+q = j} \left(  \left( \hH^{i}_p(\hC (F)) \otimes \bar{\vartriangle}(R^{\otimes (i-1)}_q) \right)
\oplus \left( Z^{i}_p(\hC(F)) \otimes R^{\otimes i}_q / \bar{\vartriangle}(R^{\otimes (i-1)}_q)  \right)  \right),
\]
where $M_p$ denotes the degree $p$ part of a graded module $M$, $\hH(\hC(F))$ denotes the homology of the restricted Bollob\'{a}s-Riordan complex and
$ Z^{i}_p(\hC (F)) = \ker \left( d^{i}(\hC (F))_{p} \right)$
are  its cycles. Moreover $H^{i+e(F)}_j ( \widetilde{\C} (F)) =  H^{i}_j ( \C (F))$.
\end{proposition}

We also note that the following deletion-contraction exact sequence holds for the homology of the restricted Bollob\'{a}s-Riordan complex.
\begin{theorem}
Let $e$ be a bridge of a fatgraph $F$ and let $F-e$ denote $F$ with the bridge $e$ deleted and $F/e$ denote $F$ contracted along the edge $e$,  then there exists a deletion-contraction exact sequence
\[
\begin{diagram}
\hH^{\ast}\left(  \hC (F - e)  \right)& & \rTo  && \hH^{\ast}\left(  \hC (F) \right)  \\
& \luTo  & &\ldTo &\\
& & \hH^{\ast}\left(  \hC (F/e) \right) \otimes V&&
\end{diagram}
.
\]
\end{theorem}

The proof of this is similar to the proof of property (5) of Proposition~\ref{basics} and is therefore excluded.

\subsection{Categorification of the multivariate Bollob\'{a}s-Riordan polynomial}
\label{sec.MBRP}
We now generalize the homology from the previous subsection to obtain a homology for the full three-variable Bollob\'{a}s-Riordan polynomial. We also prove relations between this homology theory, Khovanov homology and our chromatic homology from Section~\ref{construction}.

Notice that if in the proof of Lemma~\ref{lm:br5} we set $z= (r+r^{-1})/(q+q^{-1})^{\frac{1}{2}}$ instead of $1$ we obtain the two variable polynomial
\[\sum_{H\in \s(F)}(q+q^{-1})^{v(H)+p(H)}(-q)^{e(F)-e(H)} (r+r^{-1})^{2g(H)}.\]
We can modify it to the three variable polynomial
$$
R'(F,q,r,s) = \sum_{H\in \s(F)}(q+q^{-1})^{v(H)+p(H)}(-q)^{e(F)-e(H)} (r+r^{-1})^{2g(H)}(1+s^{-2})^{e(F)-e(H)}.
$$
It is straightforward to check that $R'$ is equivalent to the full
Bollob\'{a}s-Riordan polynomial (\ref{eq:BRpoly}).

The polynomial $R'$ can be categorified using a straight forward modification of the construction in Section~\ref{construction}.
Essentially all this involves is replacing the graded modules in the construction of the chain complex with trigraded modules.
To do this let $V$, $U$ and $R$  be the trigraded free modules of rank two with generators in degrees $(\pm 1, 0, 0)$, $(0,\pm 1, 0)$ and $(0,0,-2)$ and $(0,0,0)$ respectively.
Then we assign the module
 $V^{\otimes v(H)+p(H)}\otimes U^{\otimes 2g(H)} \otimes R^{\otimes h(H)} \{h(H)\}$ to each state. The per-edge maps are defined as before but acting in the relevant grading.
Finally, defining the graded dimension as
\[
\qdim(M)=\sum_{i,j,k} q^i r^j s^k  \, rk(M_{(i,j,k)}),
\]
we obtain the desired categorification of $R'(F,q,r,s)$.

Observe that our chromatic homology is recovered by projecting the  trigraded homology groups onto a single grading.

The idea of using multi-graded modules was also used in \cite{JHR}, to categorify the Tutte polynomial and \cite{Rong} to categorify the Bollob\'{a}s-Riordan polynomial.

\medskip

Motivated by realizations of the Jones polynomial as a signed Tutte polynomial~\cite{Ka} or as a Potts partition function~\cite{Jo} (such realizations will be discussed further in the next section),
 we extend the function
$R'$ to fatgraphs with {\em signed edges}, that is fatgraphs such that each bridge of $F$ is decorated with plus or minus sign. Let us reserve the symbol $F_s$ for  signed fatgraphs.
For a state $H$ of a signed fatgraph, let $e_-(H)$ (respectively $e_+(H)$) denote the number of bridges in $H$ with a negative (resp. positive) sign.
We  define the {\em height} of a state $H$ as $h_s(H)= e_-(F)-e_-(H)+e_+(H)$.
The construction of the homology described above with respect to this new height function, gives a categorification of the polynomial
\begin{multline*}
R'(F_s,q,r,s) =  (-q -qs^{-2})^{e_-(F_s)} \\ \sum_{H\in \s(F_s)}(q+q^{-1})^{v(H)+p(H)} (r+r^{-1})^{2g(H)}(-q (1+s^{-2}))^{e_+(H)-e_-(H) }
  .\end{multline*}
Up to normalization this can be seen to be equivalent to the more simple polynomial
\[
 \sum_{H\in \s(F)}x^{k(H)} y^{g(H)}  \prod_{e\in E(H)} z_e
, \]
where $z_e$ equals $x^{-1}z$ for an edge of positive weight and  $x^{-1}z^{-1}$ for an edge of negative weight.

\medskip

In order to write down a universal coefficient type theorem for the categorification of the three variable polynomial $R'(F_s,q,r,s)$, we introduce some notation.
Let  $HB(\D(F_s) )$ be the homology of the complex $\D(F_s)$ associated with the three variable Bollob\'{a}s-Riordan polynomial of signed graphs $R'$ described above.
Clearly $\D(F_s)= \left(  V^{\otimes v}   \otimes\hD^i(F_s) \otimes  R^{\otimes i}  , 1^{\otimes v} \otimes \hd^i  \otimes  \bar{\vartriangle}_i \right)$.

Let $\hHB$ denote the homology of the subcomplex $(\hD(F_s)^i , \hd^i)$.
Its Euler characteristic is equal to the polynomial
\[
\widehat{R'}(F_s,q,r) =  \sum_{H\in \s(F_s)}(q+q^{-1})^{p(H)} (r+r^{-1})^{2g(H)}(-q)^{h_s (H)} .
\]
\begin{proposition}\label{th:UCT}
Let $F_s$ be a  fatgraph. Then
\begin{multline*}
HB^i_{(j,k,l)}\left( \D (F_s) \right)
 = \bigoplus_{p +q = j} \left(  \left( \left( \hHB^{i}_{(p,k,0)}(\hD(F_s)) \otimes \bar{\vartriangle}(R^{\otimes (i-1)}_{(0,0,l)}) \right) 
\right. \right.  \\ \left. \left.
\oplus \left(\widehat{ Z}^{i}_{(p,k,0)}(\hD(F_s)) \otimes R^{\otimes i}_{(0,0,l)} / \bar{\vartriangle}(R^{\otimes (i-1)}_{(0,0,l)})  \right)  \right)\otimes V^{\otimes v}_{(q,0,0)}  \right),
\end{multline*}
where $M_{(j,k,l)}$ denotes the degree $(j,k,l)$ part of a trigraded module $M$ and
$  \widehat{Z}^i_{(j,k,l)}(\hD(F_s)) = \ker \left( \hd^{i}_{(j,k,l)} \right)$ are  the cycles of $\hD(F_s)$.
\end{proposition}
\begin{proof}
Let us consider the subcomplex 
 $\D'(F_s)= \left( \hD^i(F_s) \otimes  R^{\otimes i}  , \hd^i  \otimes  \bar{\vartriangle}_i \right)$.
Let $(HB')^*$ denote its homology.

By the universal coefficient theorem (see eg. \cite{Br}) 
\[
HB^i_{(j,k,l)}\left( \D (F_s) \right)
 = \bigoplus_{p +q = j} \left( (HB')^{i}_{(p,k,l)}(\D'(F_s)) \otimes V^{\otimes v}_{(q,0,0)}\right),
\]
and we see that it remains to understand $(HB')^{i}_{(p,k,l)}(\D'(F_s))$.
Let $\widehat{Z}^i = \ker (\hd^i)$ and $\widehat{B}^i = \im (\hd^{i-1})$.
For convenience we will denote $\overline{\vartriangle}_i $ by  $\overline{\vartriangle}$.
Then since  $\hd^i$ and  $\overline{\vartriangle}_i$ are of degree zero we have
\[ (HB')^{i}_{(p,k,l)}(\hD(F_s)) =
\left(
\frac{\widehat{Z}^i \otimes R^{\otimes i}}{\widehat{B}^i \otimes \overline{\vartriangle} (R^{\otimes (i-1))}}
 \right)_{ (p,k,l)}
=
\frac{ \widehat{Z}^i_{(p,k,0)} \otimes (R^{\otimes i})_{(0,0,l)} }{ \widehat{B}^i_{(p,k,0)} \otimes \overline{\vartriangle} (R^{\otimes (i-1)})_{(0,0,l)} }.
\]
Now $\widehat{Z}^i_{(p,k,0)}$ and  $\widehat{B}^i_{(p,k,0)}$ are free abelian groups with
$\widehat{Z}^i_{(p,k,0)} \subset \widehat{B}^i_{(p,k,0)}$, and so we have
 $\widehat{Z}^i_{(p,k,0)} = \bigoplus_{\alpha \in I} n_{\alpha} \Z $ and
 $\widehat{B}^i_{(p,k,0)} = \bigoplus_{\alpha \in I} m_{\alpha} \Z $, where the sum is over the same finite index $I$ and $m_{\alpha} | n_{\alpha}$ for each $\alpha \in I$. Also it is easy to show that
$\left( R^{\otimes (i-1))} \right)_{ (0,0,l)} = \Z^{\oplus N}$ and $\overline{\vartriangle} \left(R^{\otimes (i-1)} \right)_{(0,0,l)}  =\Z^{\oplus M}$ for some $M$ and $N$ ({\em i.e.} the generators of these groups are $1$).
We can then write
\[
\frac{ \widehat{Z}^i_{(p,k,0)} \otimes (R^{\otimes i})_{(0,0,l)} }{ \widehat{B}^i_{(p,k,0)} \otimes \overline{\vartriangle} (R^{\otimes (i-1)})_{(0,0,l)} }
=
\frac{\left( \bigoplus_{\alpha \in I} n_{\alpha} \Z  \right)\otimes \Z^{\oplus N}  }{\left( \bigoplus_{\alpha \in I} m_{\alpha} \Z  \right)\otimes \Z^{\oplus M}},
\]
where $M|N$.
Using standard properties of the tensor product, we see that this can be written as
\[
\bigoplus_{\alpha \in I} \frac{(n_{\alpha} \Z)^{\oplus N}   }{(m_{\alpha} \Z)^{\oplus M}  }
=
\bigoplus_{\alpha \in I} \frac{(n_{\alpha} \Z)^{\oplus(M+P)}   }{(m_{\alpha} \Z)^{\oplus M}  }
=
\left( \bigoplus_{\alpha \in I} \frac{n_{\alpha} \Z   }{m_{\alpha} \Z  }\right)^{\oplus M} \bigoplus \;
\left( \bigoplus_{\alpha \in I} n_{\alpha} \Z \right)^{\oplus P},
\]
which by definition is equal to
\begin{multline*}
\left(\widehat{Z}_{(p,k,0 )}^i / \widehat{B}_{(p,k,0)}^i\right)^{\oplus M} \oplus \left(\widehat{Z}_{(p,k,0)}^i \right)^{\oplus P}\\
  = 
 \left( \hHB^{i}_{(p,k,0)}(\hD(F_s)) \right)^{\oplus M}
\oplus \left( Z^{i}_{(p,k,0)}(\hD(F_s)) \right)^{\oplus P}
\\
 =
 \left( \hHB^{i}_{(p,k,0)}(\tD(F_s)) \otimes \Z^{\oplus M} \right)
\oplus \left( Z^{i}_{(p,k,0)}(\tD(F_s)) \otimes \Z^{\oplus P} \right).
\end{multline*} 
Now since $R^{\otimes i}_{(0,0,l)} = \Z^{\oplus N}$ and $\bar{\vartriangle}(R^{\otimes (i-1)})_{(0,0,l)} = \Z^{\oplus M} $ we have that
\[ R^{\otimes i}_{(0,0,l)} / \bar{\vartriangle}(R^{\otimes (i-1)}_{(0,0,l)}) = \Z^{\oplus (N-M)} = \Z^{\oplus P}.\] The above is equal to
\[
\left(  \left( \hHB^{i}_{(p,k,0)}(\hD(F)) \otimes \bar{\vartriangle}(R^{\otimes (i-1)}_{(0,0,l)}) \right)
\oplus \left( Z^{i}_{(p,k, 0)}(\hD(F))\otimes R^{\otimes i}_{(0,0,l)} / \bar{\vartriangle}(R^{\otimes (i-1)})_{(0,0,l)}  \right)  \right),
\]
as required.
\end{proof}

Note that in the above proof we do need to be careful about the order of the summands in the direct sums.

 Proposition~\ref{th:BRuniv} is a corollary of this result by projecting the trigraded modules onto a single graded dimension as follows.
 \begin{proof}[Proof of Proposition~\ref{th:BRuniv}.]
Regard the fatgraph $F$ as a signed fatgraph $F_s$ by assigning a negative sign to  each edge. We have
\[
\tH^i_j(\tC (F)) = \bigoplus_{p+q+r=j} HB^{i}_{(p,q,r)}( \D (F_s) ) .
\]
An application of the above proposition gives
\begin{multline*} 
\bigoplus_{p+q+r=j} \;
\bigoplus_{s+t = p} \left(  \left( \left( \hHB^{i}_{(s,q,0)}(\hD(F_s)) \otimes \bar{\vartriangle}(R^{\otimes (i-1)}_{(0,0,r)}) \right) \right. \right.\\
\oplus \left.\left.\left( Z^{i}_{(s,q,0)}(\hD(F_s)) \otimes R^{\otimes i}_{(0,0,r)} / \bar{\vartriangle}(R^{\otimes (i-1)}_{(0,0,r)})  \right)  \right)\otimes V^{\otimes v}_{(t,0,0)}  \right).
\end{multline*} 
Projection onto a single graded variable then gives 
\begin{multline*} 
\bigoplus_{p+q+r=j} \;
\bigoplus_{s+t = p} \left(  \left( \left( \hHB^{i}_{s+q}(\hD(F_s)) \otimes \bar{\vartriangle}(R^{\otimes (i-1)}_{r}) \right)\right. \right. \\
\oplus\left.\left. \left( Z^{i}_{s+q}(\hD(F_s)) \otimes R^{\otimes i}_{r} / \bar{\vartriangle}(R^{\otimes (i-1)}_{r})  \right)  \right)\otimes V^{\otimes v}_{t}  \right).
\end{multline*} 
Which after reindexing gives the required formula. 
\end{proof}

\medskip

The convenience of our homology theory for the Bollob\'{a}s-Riordan polynomial of weighted graphs is described in the following theorem. The theorem provides one of the desired Thistlethwaite-type relations between graph and knot  homologies.
\begin{theorem}
\label{thm.exx}
The Khovanov categorification of the Jones polynomial as well as the categorification of the chromatic polynomial described in the Section~\ref{construction} may be recovered from the categorification of $R^{\prime}$ just described.
\end{theorem}
\begin{proof}
We begin by using proposition~\ref{th:UCT} to determine the homology
 $\hHB(\hD(F_s))$ from $HB(\D(F_s))$.
To do this suppose that  $r_{-2}$ and $r_0$ are the generators of $R$ in graded dimension $(0,0,-2 )$ and $(0,0,0)$ respectively.

First, consider the map $f$  defined  by $r_{-2} \mapsto 0$. Applying this to the homology gives
\begin{multline*}
f\left( HB^i_{(j,k,l)}\left( \D (F_s) \right) \right) \\
 = \bigoplus_{p+q = j}
 \left(  \left( \left( \hHB^i_{(p,k,0)}(\hD(F_s)) \otimes \mathbb{Z} \right)
\oplus \left( Z^i_{(p,k,0)}(\hD(F_s)) \otimes (\mathbb{Z}/\mathbb{Z}  ) \right) \right) \otimes V^{\otimes v}_{(q,0,0)}  \right),
\end{multline*}
which can be written
\[ = \bigoplus_{p+q = j}
 \left(  \left(  \hHB^i_{(p,k,0)}(\hD(F_s))
\oplus  Z^i_{(p,k,0)}(\hD(F_s))   \right) \otimes V^{\otimes v}_{(q,0,0)}  \right),
\]
since only basis elements in graded degree $(j,k,0)$ are not killed by $f$.
Secondly, notice that
$R^{\otimes i}_{(0,0,-2i)} = \Z$ and $ \bar{\vartriangle}(V^{\otimes (i-1)})_{(0,0,-2i)} =0$, therefore  proposition~\ref{th:UCT} also gives
\[
HB^i_{(j,k,-2i)}\left( \D(F_s) \right)
 = \bigoplus_{p+q = j}  Z^i_{(p,k,0)}(\hD(F_s)) \otimes \Z   \otimes V^{\otimes v}_{(q,0,0)}.
\]
Finally, since we know each free module  $V^{\otimes v}_{(q,0,0)} $ we know $Z^i_{(p,k,0)}(\hD(F_s))$ and hence $ \hHB^i_{(p,k,0)}(\hD(F_s))$.

The Khovanov homology of an associated link can be recovered from this as in the
 discussion in the discussion that will follow in Section~\ref{ss:khovanov}: associate the  crossing
\rotatebox{90}{\includegraphics[height=3mm]{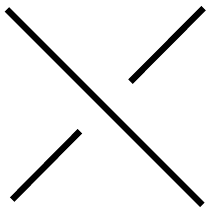}}
to a edge of negative weight and
\includegraphics[height=3mm]{tosmooth} 
to an edge of positive weight.
The height then equals the number of 1-smoothings of the associated link. The Khovanov homology is then a normalization of $\hHB$ after projection onto a single grading.

The second statement follows by regarding an unsigned  fatgraph as a signed fatgraph whose edges all have negative weight and projecting onto a single graded dimension. 
\end{proof}

\section{Independence of planar embeddings}\label{sec:khov}

We begin this section  by  describing the relationship of our homology theory with Khovanov's categorification of the Jones polynomial \cite{Kh}. This provides another Thistlethwaite-type relation between graph and knot homologies. We will
describe this relation in Subsection~\ref{ss:khovanov} and apply it in Subsection~\ref{ss:proof} the prove the following theorem on the independence of our chromatic homology of a plane graph on the choice of planar embedding.
\begin{theorem} \label{th:invariance}
Let $F$ and $F^{\prime}$ be two genus $0$ fatgraphs with the same associated  graph $G$, and let $P$ be the Poincar\'{e} polynomial of the homology. Then
\[ P(F)=P(F^{\prime}),\]
{\em i.e.} the Poincar\'{e} polynomial is independent of the  embedding of the graph $G$.
\end{theorem}

\subsection{The Relation to Khovanov Homology and Knots}
\label{ss:khovanov}

Let $F$ be a genus $g$ fatgraph. As mentioned earlier, $F$ is equivalent to  a genus $g$ surface which we will denote $\Sigma_g$.  $F$ gives rise to an alternating link $L\subset \Sigma_g \times I$, and a canonical diagram onto $F$ by associating a crossing
\includegraphics[height=3mm]{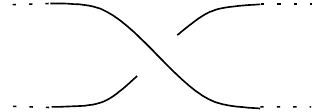} 
to each bridge
\includegraphics[height=3mm]{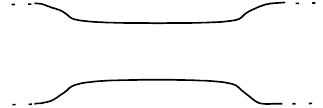} 
and connecting these crossings according to the cyclic ordering at the islands of the fatgraph.
We call this the {\em associated link}.
The following figure shows a fatgraph with one island of degree 4 and one island of degree 2, and  its associated link.
\[ \includegraphics[width=4cm]{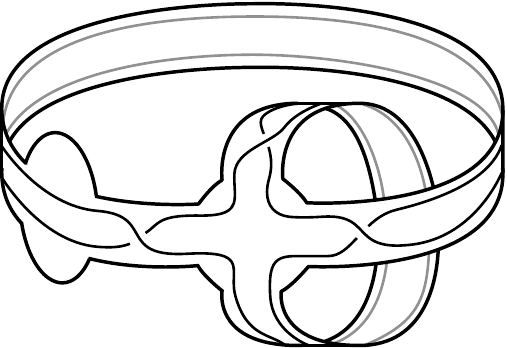} 
\]

Just as with link diagrams on $S^2$, we can consider the  {\em smoothing} of a crossing. A {\em 0-smoothing} is defined locally on a link diagram by changing a crossing which looks like 
\includegraphics[height=3mm]{tosmooth} 
 to look like
 \includegraphics[height=3mm]{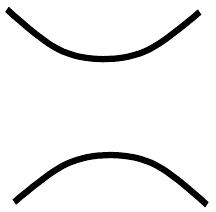}; 
  and for a {\em 1-smoothing} replacing the crossing with
  \includegraphics[height=3mm]{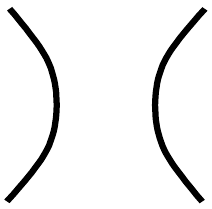} 
   \,. A {\em state} of a link diagram  is what is obtained by smoothing all of the crossings of the link diagram.
There is a clear correspondence between the states of a fatgraph and the states of the associated link. This is summarized in the  table below.
\begin{center}\begin{tabular}{|c|c|c|}
\hline
  smoothing &  ass. link &  fatgraph   \\
\hline
\raisebox{11mm}{} \raisebox{4mm}{0} &
\includegraphics[height=1cm]{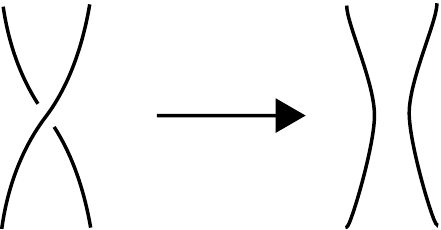}
  & 
 \includegraphics[height=1cm]{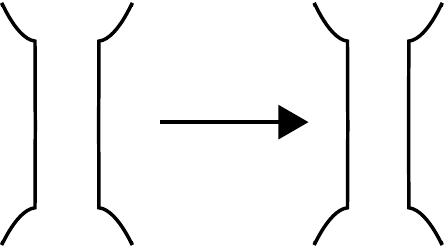}
  \\
\hline
\raisebox{11mm}{} \raisebox{4mm}{1}&
 \includegraphics[height=1cm]{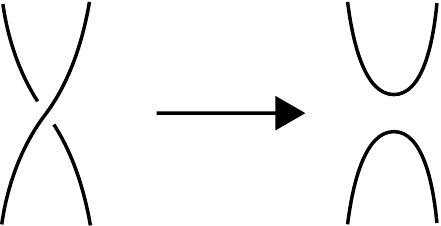}
 & \includegraphics[height=1cm]{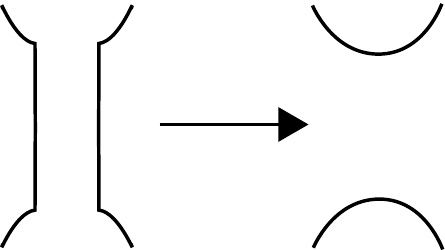}
 \\
\hline
\end{tabular}\end{center}

\medskip

From now on we will restrict ourselves to genus $0$ fatgraphs. The associated links can then be regarded as links in $S^3$.

The Khovanov homology of a link is constructed in essentially the same way as the complex in Section~\ref{construction}: given a link diagram $L$ we define the {\em height}, $h(S)$, of a state $S$ to be the number of 1-smoothings used in its construction and we let $p(S)$ denote the number of cycles in the state $S$. The chain modules are constructed by assigning the graded module $V^{\otimes p(S)} \{ h(S) \}$, which has graded dimension $q^{h(S)}(q+q^{-1})^{p(S)}$. Again the $i$-th chain group is defined to be the direct sum of all modules assigned to states of height $i$.  Just as in Section~\ref{construction}, one can move from a state of height $i$ to $i+1$ by merging or splitting cycles and we use the per-edge maps $\widetilde{m}$ and $\widetilde{\vartriangle}$ from Section~\ref{construction} to move between the corresponding modules. The differentials are then obtained by summing over all of the per-edge maps of the appropriate height as before. 
This gives a chain complex $ \tC (L) = \left(  \tC^i(L) , \td^i  \right)$. We let $\tHK(L)$ denote
the  homology of this complex.

If $n_{\pm}$ denotes the number of $\pm$-crossings of $L$ (the sign of a crossing will be defined in Subsection~\ref{ss:proof}), then the Khovanov complex is defined as the following normalization of $\tC$:
\[ \C (L) = \left(  \C^i(L) , d^i  \right)  = \tC (L)[-n_- ]\{ n_+ - 2n_- \} = \left(  \tC^i(L) , \td^i  \right)[-n_- ]\{ n_+ - 2n_- \} . \]
The homology of this complex is called {\em Khovanov homology}, $HK(L)$. It can be shown (\cite{Kh, BN1}) that the Euler characteristic $\chi (HK(L))$ is equal to the Jones polynomial of $L$ and that the homology itself is a knot invariant which is strictly stronger than the Jones polynomial.

observe that by the correspondence between the states of a fatgraph and the states of an associated link described above, we have \[\tC (F) = \left(  V^{\otimes v}   \otimes\tC^i(L) \otimes  R^{\otimes h}  , 1^{\otimes v} \otimes \td^i  \otimes  \bar{\vartriangle}_h \right),\] where $L$ is the associated link of a planar fatgraph.

This observation leads to the following Thistlethwaite-type theorem, which is  a corollary of Theorem~\ref{th:BRuniv}.
\begin{theorem} \label{th:univcoeff}
Let $F$ be a genus $0$ fatgraph and $L \subset S^3$ be the associated link with an arbitrary orientation. Then
\begin{multline*}
H^i_j\left( \widetilde{\C}(F) \right) \\
 = \bigoplus_{p+q+r = j} \left( \left( \left( \widetilde{HK}^i_p(L) \otimes \bar{\vartriangle} R^{\otimes (i-1)}_q \right)
\oplus \left( Z^i_p(L) \otimes R^{\otimes i}_q / \bar{\vartriangle} R^{\otimes (i-1)}_q  \right)\right) \otimes V^{\otimes v}_r  \right),
\end{multline*}
where $M_p$ denotes the degree $p$ part of a graded module $M$,
\[ Z^i_p(L) =   \ker(\td^i_p)
\;\; \left( = \ker \left( d_{p + 2n_{-} - n_+}^{i+n_{-}} \right)\right) \] 
are  cycles determined by the Khovanov complex and for the planar fatgraph $F_s$ corresponding to $L$,
\[  \widetilde{HK}^i_p(L) = HK_{p + 2n_{-}- n_+}^{i+n_{-}} (L),\]
where $HK$ denotes Khovanov homology.
Moreover $H^i_j(\C)= H^{i+e(F)}_j(\tC)$.
\end{theorem}
Note that one can prove an analogous result for the full ({\em i.e.} without the restriction to alternating links) Jones polynomial using the homology of the signed Bollob\'{a}s-Riordan polynomial.

\subsection{Proof of Theorem~\ref{th:invariance}} \label{ss:proof}
This subsection is devoted to the proof of Theorem~\ref{th:invariance}.
Our method is to reduce the graph theoretical  problem to one of knot theory and to prove the result using this relation.

First we need to understand how the two fatgraphs $F$ and $F^{\prime}$ and their corresponding links are related. For  this we find it convenient to switch to the language of embedded graphs. Recall that a genus 0 fatgraph is equivalent to an embedding of the associated graph $G \hookrightarrow S^2$.
We will need the following two local moves on embedded graphs.
Let $G \subset S^2$ be a connected embedded graph.
A {\em 1-flip} is a move which replaces a 1-connected component of the map $G$ with its rotation by $\pi$ around the axis in the $xy$-plane which intersects the 1-connecting vertex. For example,
\[ \includegraphics[height=2cm]{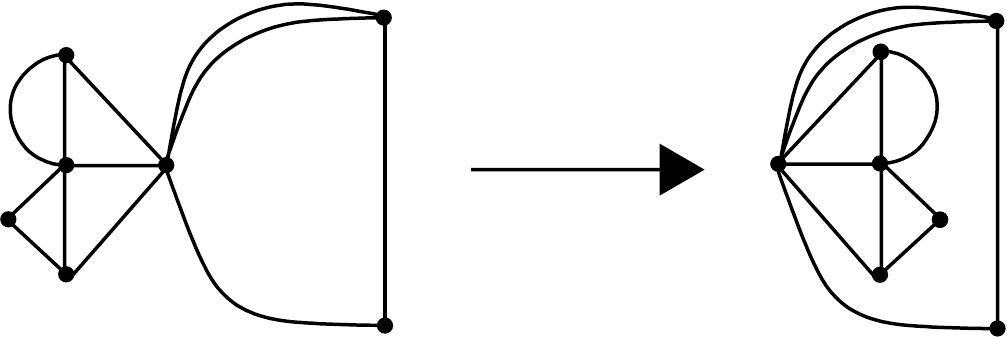}
.\]

A {\em 2-flip} is a move which replaces a 2-connected component of $G$ with its rotation by $\pi$ around the axis determined by its 2-connecting vertices. For example,
\[ \includegraphics[height=2cm]{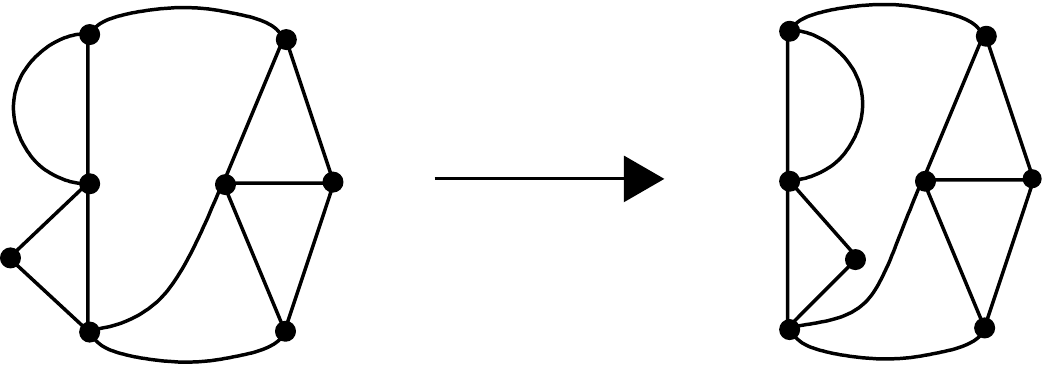}
.\]

The following theorem relates two planar embeddings of the same graph (see \cite{MM}, or \cite{MT}).

\begin{theorem}
\label{thm.moh}
Let $G$ be a  connected graph and $f, \, f^{\prime} : G \rightarrow S^2$ be two planar embeddings. Then $f(G)$ and $f^{\prime}(G)$ are related by a sequence of 1-flips and 2-flips.
\end{theorem}

We need to understand how a flip changes the associated link.
First consider the 1-flip. By regarding a 1-connected graph  as two components with a vertex identified, it is easy to see that the associated link is non-prime with a connect sum determined by the connecting vertex, and that the embeddings of the two components determine where the connect sum occurs. (Recall that the sum of two links is the link formed by cutting open an arc of each link and identifying the free ends in a way consistent with orientation. Although this process is well defined with respect to isotopy for knots, in general it will depend upon which components of the link have been identified under the sum.)

It is then clear that in terms of the associated link, a 1-flip simply changes the way we connect the two links in the connect sum. We will discuss this in more detail in the proof of Lemma~\ref{lem:linkcon}.

\medskip

Now, by considering the example
\[  \includegraphics[height=4cm]{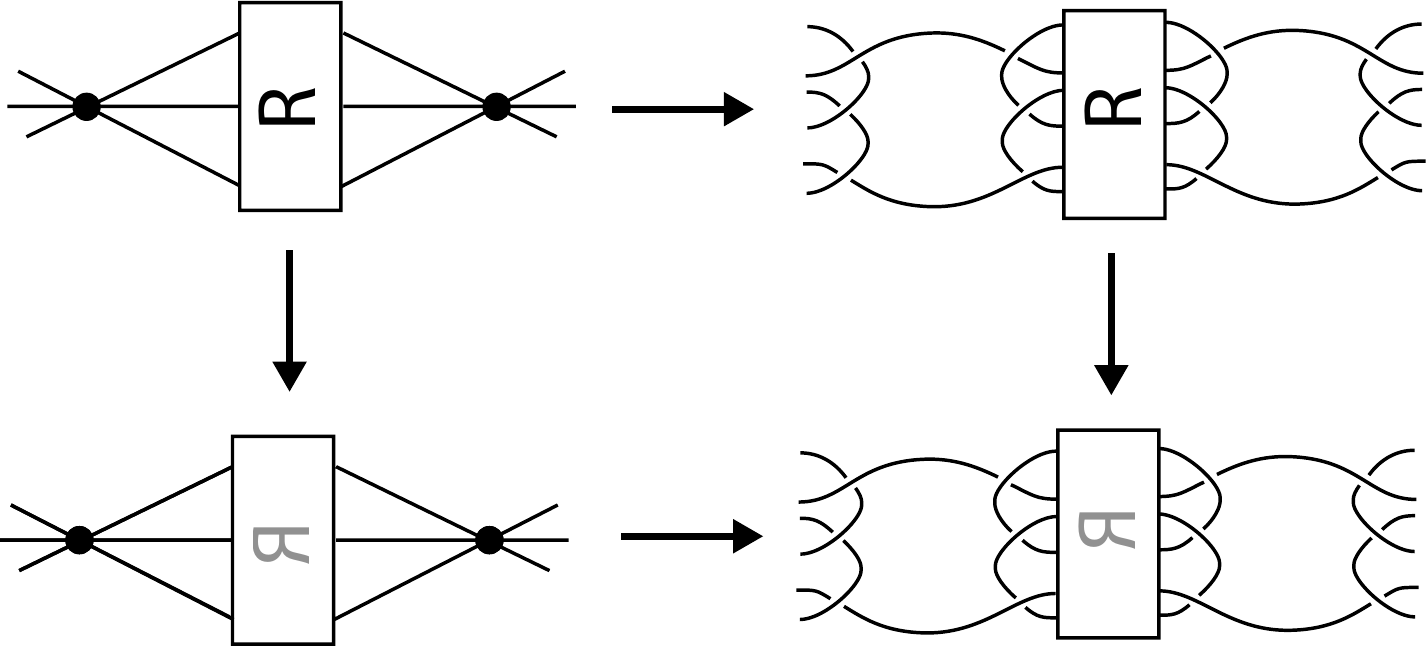}
,\]
it is easy to see that a 2-flip corresponds to a move which replaces a $2$-tangle with its rotation by $\pi$ around the axis between the tangle ends:
\[
 \includegraphics[height=3cm]{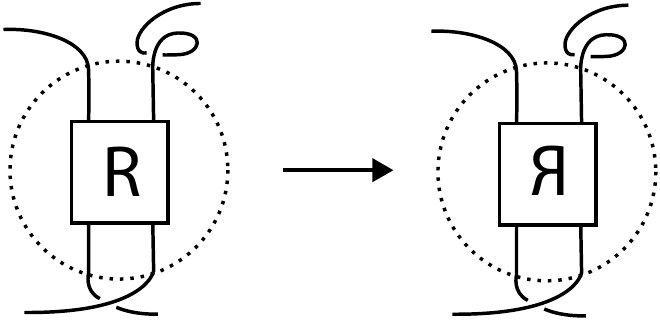}
.
\]
This is a form of Conway mutation of the link (\cite{Co}). We  refer to this specific move on the link diagram simply as {\em mutation}.

This discussion gives the following lemma.
\begin{lemma} \label{lem:relation}
Two genus 0 fatgraphs have the same associated connected graph if and only if their associated links are obtained as the connect sum of the same set of links and by a  sequence of mutations.
\end{lemma}

So far we have only discussed un-oriented links, but our application of knot theory requires a choice of orientation of the links. We  need to be careful in this choice of orientation.
Recall that the \emph{sign} of a crossing of an oriented link is the assignment of $\pm 1$  according to the following scheme:
\[
 \includegraphics[height=1.6cm]{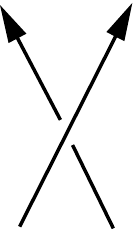}
 \; \raisebox{6mm}{+1 \; ,}
\hspace{2cm}
\includegraphics[height=1.6cm]{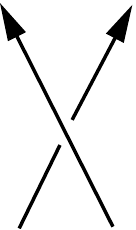}
 \; \raisebox{6mm}{-1}.
\]
If the components of a link $L$ are labelled $\{ 1, \ldots , n \}$, then the {\em linking number} $lk_L(i,j)$ is defined to be the sum of the signs over all crossings between the components labelled $i$ and $j$.

Given two links which are related as in Lemma~\ref{lem:relation}, we require that the corresponding crossings in the two links have the same sign.
To do this we orient  the summands of links arbitrarily.
First we deal with the case of a connected sum. 
If a connect sum requires the reversal of the orientation of a component, then we reverse all of the components of that summand. The case for mutation is a little more complicated.
 Let $R$ denote the tangle on which the mutation acts. Then, if the two free arcs at the top $R$ are both oriented into or out of the tangle, we retain all orientations. If one of the  arcs at the top $R$ is oriented into the tangle and one out of the tangle, then we reverse the orientations of {\em all} the components of the tangle $R$. See  the figure below.
 We call such an orientation of the mutant the {\em induced orientation}. It is immediate that the sign of each crossing before and after the mutation is the same.
\[
\begin{array}{ccc}
\includegraphics[width=10mm]{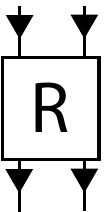}
& &\includegraphics[width=10mm]{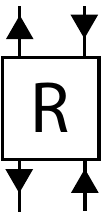}
 \\
 \text{{\em keep orientation.}} & & \text{{\em change orientation.}}
\end{array}
\]
Further to this, the labeling of $R$ induces a labeling of the components after the mutation which we call  the {\em induced labeling}.
We can now prove the theorem.
\begin{proof}[Proof of Theorem~\ref{th:invariance}.]
Let $f(G)$ and $f^{\prime}(G)$ be the embeddings corresponding to  the associated graphs $F$ and $F^{\prime}$.
The following lemma reduces the problem.
\begin{lemma}
Let $L$ and $L^{\prime}$ be the two links associated with the fatgraphs  $F$ and $F^{\prime}$. Then to prove the theorem, it is enough to show that $P\left( HK(L) \right) = P\left( HK(L^{\prime}) \right)$.
\end{lemma}

\begin{proof}
Assume $P\left( HK(L) \right) = P\left( HK(L^{\prime}) \right)$. 
 Recall that if $H = \left( H^i \right)_{i \in \mathbb{Z} }$ is the homology of some chain complex of graded $\mathbb{Z}$-modules, the Poincar\'{e} polynomial is defined by
\[
P(H) = \sum_{i \in \mathbb{Z}} t^i \qdim\left( H^i \right).
\]
The Poincar\'{e} polynomial encodes all of the torsion-free information of the homology groups.
By Theorem~\ref{th:univcoeff}, to prove the lemma it suffices to show that $\rk \left( \tHK(L)^i_j \right)
= \rk \left( \tHK(L')^i_j \right)$, and that $\rk \left( \ker (\td^i_j)  \right)= 
\rk \left( \ker (\td^i_j)  \right)$.

Clearly, in any graded degree we have
\[\rk \left(  HK^i_j \right) =
\rk \left( \ker (d^i_j)  \right)
 -
\rk \left( \im (d^{i-1}_j) \right),
\]
and by classic linear algebra we know that
\[
\rk \left( \C^i_j \right) =
\rk \left( \ker (d^i_j)  \right)
 +
\rk \left( \im (d^{i}_j) \right).
\]

Suppose we know the ranks of $\C^k_j$ and $HK^k_j$ for some $k$. If in addition to this we know the value $\rk \left( \im (d^{k-1}_j) \right)$, then by the above we can determine the values of $\rk \left( \ker (d^k_j)  \right)$ and $\rk \left( \im (d^{k}_j) \right)$.
Therefore if we know the ranks of each $\C^i_j$ and $HK^i_j$ and one value of $\rk \left( \im (d^{k-1}_j) \right)$ then we can determine every value $\rk \left( \ker (d^i_j)  \right)$ and $\rk \left( \im (d^{i}_j) \right)$ and we then know the torsion-free information of the entire complex.

It is easily seen (recall  that $v+p+2g-e= 2k$) that each term $\rk \left( \C^i_j \right)$ is equal for the two chain complexes from the two associated links $L$ and $L^{\prime}$.  
Also by assumption we have that each value of  $\rk \left(  HK^i_j \right)$ is the same in the two complexes. The above argument then tells us that if there exists a value $\rk \left( \im (d^{k_l}_l) \right)$ which is equal for the two complexes associated with $L$ and $L^{\prime}$ for each graded dimension $l$, then the torsion-free parts of the two complexes are equal.
Finally since the chain complexes in any graded degree only have finitely many non-zero terms, clearly such values  $\rk \left( \im (d^{k_l}_l) \right) =0 $ exist, completing the proof of the lemma.
\end{proof}

\medskip

We need to use the following theorem of Lee, which proved a  conjecture of Garoufalidis from \cite{Ga}.
\begin{theorem}[\cite{Le2}]
  For a reduced alternating link, $P(HK)$ is determined by the Jones polynomial, the signature of the link and the linking numbers.
\end{theorem}
Clearly the associated links are alternating. Since we are only interested in determining $P(HK)$ and Khovanov homology is a link invariant, we can reduce the associated link (recall that a link projection is said to be {\em reduced} if four distinct regions meet at every crossing).
It is well known and easily seen that Conway mutation does not change the Jones polynomial or the signature of the link (for a definition of the signature of a link see, for example, \cite{Li}) and since the Jones polynomial is multiplicative and the signature additive under the connect sum, we see that these two invariants are equal for both our associated links.
However, the linking numbers do change under the operations.
By Theorems~1.2 and 4.5 of \cite{Le1} and Corollary~A.2 of \cite{Ga}, we see that it is enough to show that if the components of each of the associated links are labelled $1, \ldots , n$ then the following formula is equal for both links,
\begin{equation}\label{eq:lees}
\sum_{ E \subset \{ 2, \ldots, n  \}}   \gamma^{ \sum_{j \in E, k \not \in E  } lk_{jk}}.\end{equation}
The following two lemmas will complete the proof of Theorem~\ref{th:invariance}.
\begin{lemma}\label{lem:linkcon}
Equation~\ref{eq:lees} does not depend upon which two components the connect sum operation acts.
\end{lemma}
\begin{proof}
Suppose we have two links $L$, with components labeled $a_1 , \ldots , a_n$, and $L^{\prime}$ with components labeled $b_1 , \ldots , b_m$.
Let $A$ be the link obtained by connect summing with respect to the components $a_1$ and $b_1$ labeling this new component $\alpha$ and $B$ be the link obtained by connect summing with respect to the components $a_2$ and $b_1$ labelling this new component $\beta$.
If $n \leq 2$ and $m=1$ the result is obvious, so assume that this is not the case.
Now if
$E \subset \{ \alpha , a_2, \ldots a_n,  b_2, \ldots b_m  \}$
then construct a subset $F$ from $E$ by replacing the element $\alpha$ with $a_1$ and $a_2$ by $\beta$, if $\alpha$ or $a_2$ are in $E$.
Similarly if
$E \subset \{ \beta , a_1, a_3, \ldots a_n,  b_2, \ldots b_m  \}$
then construct a subset $F$ from $E$ by replacing the element $\beta$ with $a_2$ and $a_1$ by $\alpha$, if $\beta$ or $a_1$ are in $E$.
Then since
\[ lk_{A}(\alpha, a_i) = lk_{L}(a_1, a_i), \quad
 lk_{A}(\alpha, b_i) = lk_{L^{\prime}}(b_1, b_i), \]\[
lk_{A}(a_i, a_j) = lk_{L}(a_i, a_j), \quad
lk_{A}(b_i, b_j) = lk_{L^{\prime}}(b_i, b_j)\]
and
\[ lk_{B}(\beta, a_i) = lk_{L}(a_2, a_i), \quad
 lk_{B}(\beta, b_i) = lk_{L^{\prime}}(b_1, b_i), \]\[
lk_{B}(a_i, a_j) = lk_{L}(a_i, a_j), \quad
lk_{B}(b_i, b_j) = lk_{L^{\prime}}(b_i, b_j)\]
and all other linking numbers are zero,
we have that
\[ \gamma^{ \sum_{j \in E, k \not \in E  } lk_{jk}}
= \gamma^{ \sum_{j \in F, k \not \in F  } lk_{jk}},
\]
and therefore for each summand in the equation~\ref{eq:lees} for the link $A$, there is a corresponding summand of equal value in the equation for the link $B$ and vice-versa.
This completes the proof of Lemma~\ref{lem:linkcon}. 
\end{proof}

Note that by our choice of orientation of the summands, to prove the theorem we are allowed to assume that the connect sum is consistent with the orientation.

\begin{lemma}
If $L$ is a link, $L^{\prime}$ is obtained from $L$ by a mutation and $L^{\prime}$ has the canonical orientation and labeling, then the sum of formula~\ref{eq:lees} is equal for the two links $L$ and $L^{\prime}$,
\end{lemma}
\begin{proof} 
We can regard the link $L$ as the identification of two $2$-tangles, $R$ and $T$. We may assume that $L^{\prime}$ is obtained from $L$ by a mutation which `flips over' the tangle $R$. Then the two links $L$ and $L^{\prime}$ and their linking numbers differ according to how the two tangles $R$ and $T$ are joined. This can be represented by the following figure.
\[
\begin{array}{ccc}
\includegraphics[width=3.5cm]{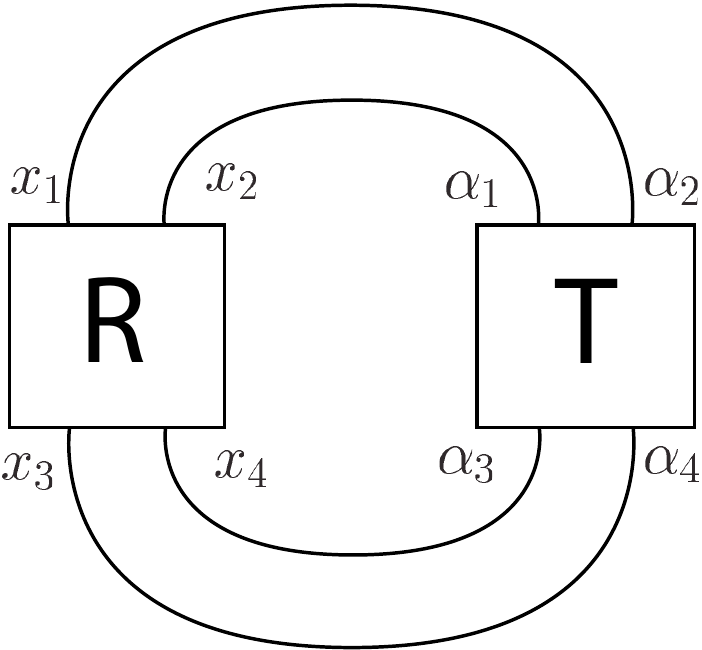}
& \hspace{1cm} & \includegraphics[width=3.5cm]{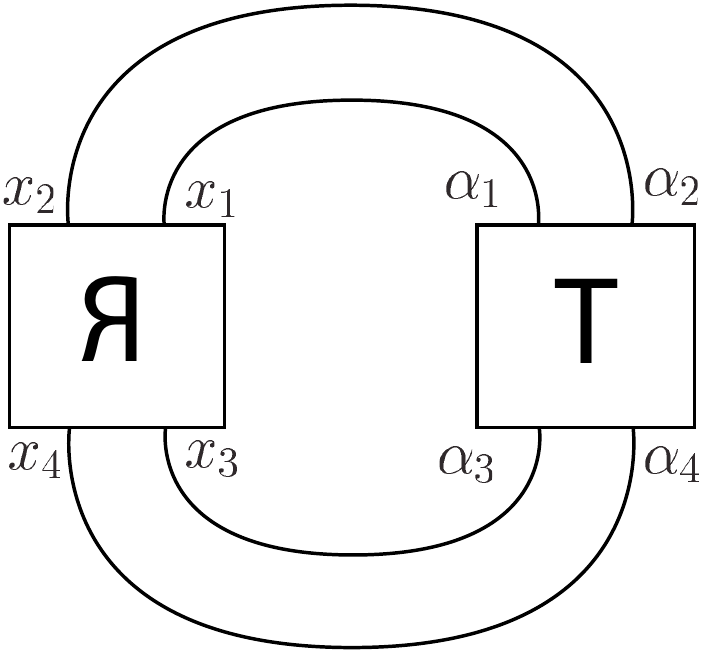}
\end{array}
\]
The eight free ends of the tangles belong to one or two components of the final link. The remainder of the proof is split into several cases according to which components the free ends will belong to in the corresponding link. Note that the number of components of $L$ and $L^{\prime}$ will always be equal.

Now suppose that the free ends of the tangle are labelled as in the above figure.
First note that only the linking numbers which involve a component coming from one of the free ends can change under mutation.

\textbf{Case 1.} Suppose that after identification all of the free ends belong to one component labelled $A$.
Then
\[
lk_{L} (A,i) =
 lk_{L^{\prime}} (A,i)
,  \]
for an arbitrary component $i$. Therefore Lee's formula is unchanged by mutation.

\textbf{Case 2.} Suppose that after identification $x_1$, $x_2$, $\alpha_1$ and $\alpha_2$ belong to the same component labelled $A$ and $x_3$, $x_4$, $\alpha_3$ and $\alpha_4$ belong to the same component labelled $B$. Then

\[
\begin{array}{ll}
lk_{L} (A,i) =& lk_{L^{\prime}} (A,i)
\\
lk_{L} (B,i) = & lk_{L^{\prime}} (B,i)
\end{array}
 \]
for an arbitrary component $i$, and again Lee's formula is unchanged by mutation.

\textbf{Case 3.} Suppose that after identification $x_1$, $\alpha_2$, $x_3$ and $\alpha_4$  belong to the same component
labelled $A$ and  $x_2$, $\alpha_1$, $x_4$ and $\alpha_3$ belong to the same component labelled $B$.
This case is more complicated. We start with some notation. Let $A_R$ be the segment of $A$ in $R$ with open ends $x_1,x_3$
and let $A_T,B_R,B_T$ be defined analogously. Hence the components $A_RA_T, B_RB_T$ of $L$ are transformed into
components $A_RB_T, B_RA_T$ of $L'$. If $n=2$ then the components of $L$ are exactly $A,B$ and Lemma is simply true.
Hence assume further that $n>2$ and $l_1=C$ is a component of $L$ different from both $A,B$. Without loss of generality
assume $C$ is a component of $R$. Let us denote by ${\mathcal L}_R$ (${\mathcal L}_T$)the set of the components of $L$ in $R$
($T$ respectively).

 Let ${\mathcal E}(L)$ be the set of all $E\subset \{2,\dots, n\}$ such that $A\in E$ and $B\notin E$,
or $A\notin E$ and $B\in E$. Analogously let ${\mathcal E}(L')$ be the set of all $E\subset \{2,\dots, n\}$ such that
$A_RB_T\in E$ and $B_RA_T\notin E$, or $A_RB_T\notin E$ and $B_RA_T\in E$.

For $E\in {\mathcal E}(L)$ let
\[
l(E)= \sum_{i\in \{A,B\}\cap E, \, j\notin E\cup\{A,B\}} lk_L(i,j) +
\sum_{i\in E-\{A,B\}, \, j\in \{A,B\}-E} lk_L(i,j).
\]

Analogously for $E\in {\mathcal E}(L')$ let
\begin{multline*}
l(E)= \sum_{i\in \{A_RB_T, B_RA_T\}\cap E, \, j\notin E\cup\{A_RB_T, B_RA_T\}} lk_{L'}(i,j) \\+
\sum_{i\in E-\{A_RB_T, B_RA_T\}, \, j\in \{A_RB_T, B_RA_T\}-E} lk_{L'}(i,j).
\end{multline*}

Clearly, it suffices to show:

{\bf Claim~1.}
\[
\{l(E);E\in {\mathcal E}(L)\}= \{l(E);E\in {\mathcal E}(L')\}.
\]

For each $X\subset \mathcal {L}_R$ containing $C$ and for each $Y\subset \mathcal {L}_T$ let
\[\mathcal{E}(L,A_R,X,Y)= \{ E\in \mathcal{E}(L); E=\{A_RA_T\}\cup X\cup Y\}\] and
\[\mathcal{E}(L',A_R,X,Y)= \{ E\in \mathcal{E}(L'); E=\{A_RB_T\}\cup X\cup Y\}.\]
We define $\mathcal{E}(L,B_R,X,Y), \mathcal{E}(L',B_R,X,Y)$ analogously.

Clearly the sets $\mathcal{E}(L,A_R,X,Y), \mathcal{E}(L,B_R,X,Y)$,
$X\subset \mathcal {L}_R$ containing $C$ and $Y\subset \mathcal {L}_T$ form a partition of $\mathcal{E}(L)$.
An analogous statement holds also for $L'$. Hence Claim~1 follows from the following.

{\bf Claim~2.} For each $X\subset \mathcal {L}_R$ containing $C$ and $Y\subset \mathcal {L}_T$,
\[
\{l(E);E\in \mathcal{E}(L,A_R,X,Y)\}=\{l(E);E\in \mathcal{E}(L,A_R,X,\mathcal{L}_T-Y)\}
\]
and the same is true when $A_R$ is replaced by $B_R$.

Claim~2 may be verified by checking. In fact, both sides are equal to
\[
\{lk_L(A_R,\mathcal{L}_R-X)+lk_L(B_r,X), lk_L(A_T, \mathcal{L}_T-Y)+ lk_L(B_T, Y)\}.
\]

This finishes the proof of Claim~2, and Claim~1 and  therefore Case~3.

\textbf{Case 4.} The final case is when,  after identification, $x_1$, $\alpha_2$, $x_4$ and $\alpha_3$  belong to the same component
labelled $A$ and  $x_2$, $\alpha_1$, $x_3$ and $\alpha_4$ belong to the same component labelled $B$.
This case is similar to Case~3 and the proof is omitted.

This finishes the proof of the lemma. 
\end{proof} 
This completes the proof of Theorem~\ref{th:invariance}. 
\end{proof} 

\begin{remark}
With some easy changes, the proof of Theorem~\ref{th:invariance} above shows that the torsion-free part of the  homology $HB(\D(F_s))$ is also invariant of the choice of genus $0$ embedding of a signed fatgraph all of whose edges are of the same sign.
\end{remark}

\section{Helme-Guizon and Rong's chromatic homology}
In this final section we provide a categorification of the Bollob\'{a}s-Riordan polynomial which unites Helme-Guizon and Rong's categorification for the chromatic polynomial (\cite{HGR1}) and Khovanov homology. As mentioned in the introduction, this addresses a question posed in \cite{HGR1}. We begin by recalling Helme-Guizon and Rong's categorification of the chromatic polynomial.

Helme-Guizon and Rong categorify the chromatic polynomial \[M(G, 1+r) = \sum_{H\in \s(G)} (-1)^{e(H)}(1+r)^{k(H)}.\]
This polynomial is categorified by considering modules $M$ which are free, graded, rank $2$, $\Z$-modules with generators $m_0$ in graded degree zero and $m_1$ in degree 1,
so that $\qdim = 1+r$. The height function is $|H|$ and the the module $M^{\otimes k(H)}$ is attached to each state. The per-edge maps are either the identity or the map induced by the  degree zero multiplication
 $m'(m_0, m_0) = m_0$, acting on merging connected components.

\medskip

Consider again the Bollob\'{a}s-Riordan  polynomial  
\[
R(F,x,y,z)  = \sum_{H}x^{r(F)-r(H)}y^{n(H)}z^{k(H)-p(H)+n(H)}.
\]
Using the definitions of the rank and nullity we can write this as
\[
 x^{-k} (yz)^{-v}
\sum_{H}(xyz^2)^{k(H)}(yz)^{e(H)}z^{-p(H)}.
\]
Setting $x=(1+r)(-1-q^{-2})$, $y=-1-q^2$, $z=(q+q^{-1})^{-1}$ and forgetting about the normalization, we obtain the polynomial
\[
B(F,q,r)=\sum_{H}(1+r)^{k(H)}(-q)^{e(H)}(q+q^{-1})^{p(H)}.
\]
We concern ourselves with this evaluation of the Bollob\'{a}s-Riordan polynomial, which we note is not equivalent to the $3$-variable Bollob\'{a}s-Riordan polynomial

The polynomial $B$ can be categorified. The chain complex is constructed using rank 2, free, bigraded modules $V$ and $M$ with basis generators $v_{\pm}$ in graded degree $(\pm 1, 0)$ and generators $m_0$ and $m_1$ in graded degrees $(0,0)$ and $(0,1)$ respectively. The module $M^{\otimes k(H)}\otimes V^{\otimes p(H)}$ is assigned to each state and the per-edge maps are of the form 
$m' \otimes \widetilde{m}$  and  $m' \otimes \widetilde{\vartriangle}$, where $\widetilde{m}$, $ \widetilde{\vartriangle}$ and $m'$ are the obvious bigraded versions of the maps defined in Section~\ref{construction} and above.

Notice that for a planar fatgraph $F$, this chain complex is of the form $E(F)=\left( C^i \otimes D^i, d^i \otimes \partial^i\right)$ where $C(F)=\left( C^i, d^i\right)$ is the chain complex of Helme-Guizon and Rong's chromatic homology and, by Subsection~\ref{ss:khovanov}, $D(F)=\left( D^i,  \partial^i\right)$ is the unnormalized Khovanov complex of the reflection (since the height function here corresponds to the addition of edges rather than the removal) of the associated link. 
Also note that the maps $f:M\rightarrow 1$ and $g: V \rightarrow 1$ clearly induce chain maps $f: E(G) \rightarrow D(G)$ and  $g: E(G) \rightarrow C(G)$. Putting all of this together we obtain:
\begin{proposition}
For a planar fatgraph $F$ there is a homology theory for the Bollob\'{a}s-Riordan polynomial which comes equipped with two natural homomorphisms,  one to the chromatic homology of Helme-Guizon and Rong and the other to the Khovanov homology of the reflection of the associated link.
\end{proposition}

\appendix

\section{Table of  polynomials}
To simplify the entries of the table we write $v,p,\ldots$ for $v(H), p(H), \ldots$ and sums are taken over the appropriate set of states.
\begin{center}\begin{tabular}{|c|c|c|c|}
\hline
 polynomial &  state sum &  complex & homology   \\
\hline
 $Z(F,q)$ &  $ (-1)^{e(F)}\sum(q+q^{-1})^{v+p+2g}((-q)(1+q^{-2}))^{h}$ &  $\C(F)$ & $H(\C(F))$   \\
\hline
 - &  $\sum(q+q^{-1})^{v+p+2g}((-q)(1+q^{-2}))^{h}$ &  $\tC(F)$ & $\tH(\tC(F))$   \\
\hline
 $\widehat{R}(F,q)$ &  $ \sum(q+q^{-1})^{v+p+ 2g}(-q)^{h}$ &  $\hC (F)$ & $\hH(\hC (F))$   \\
\hline
 $R'(F_s,q,r,s)$ &  $ \sum(q+q^{-1})^{v+p} (r+r^{-1})^{2g}(-q (1+s^{-2}))^{h_s }$ &  $\D (F_s)$ & $HB(\D (F_s))$   \\
\hline
 $\widehat{R'}(F_s,q,r)$ &  $\sum(q+q^{-1})^{p} (r+r^{-1})^{2g}(-q)^{h_s }$ &  $\hD (F_s)$ & $\hHB(\hD (F_s))$   \\
\hline
 Jones poly.  &  $ (-1)^{n_-}q^{n_+ - 2n_-} \sum(q+q^{-1})^{p}(-q)^{h}$ &  $\C (L)$ & $HK(L)$  \\
\hline
 - &  $\sum(q+q^{-1})^{p}(-q)^{h}$ &  $\tC (L)$ & $\tHK(L)$  \\
\hline
\end{tabular}\end{center}

\end{document}